\documentclass{article}

\usepackage{authblk}
\usepackage[numbers]{natbib}
\usepackage{geometry}
\usepackage{amsmath,amssymb}

\usepackage{orcidlink}

\usepackage{multicol}

\usepackage{graphicx}
\usepackage{epstopdf}
\usepackage[ruled,vlined,onelanguage,linesnumbered]{algorithm2e}
\usepackage{enumitem}
\PassOptionsToPackage{hyphens}{url}\usepackage{hyperref}

\usepackage{cleveref}
\crefformat{equation}{(#2#1#3)}

\newcommand{\R}{\mathbb{R}} 

\newcommand{\cB}{\mathbf{\mathcal{B}}}
\newcommand{\cC}{\mathbf{\mathcal{C}}}
\newcommand{\cD}{\mathbf{\mathcal{D}}}
\newcommand{\cS}{\mathbf{\mathcal{S}}}
\newcommand{\cT}{\mathbf{\mathcal{T}}}
\newcommand{\cX}{\mathbf{\mathcal{X}}}
\newcommand{\cY}{\mathbf{\mathcal{Y}}}

\newcommand{\vx}{\mathbf{x}}
\newcommand{\vy}{\mathbf{y}}

\newcommand{\vR}{\mathbf{R}}
\newcommand{\vEin}{\mathbf{E}_{\text{in}}}
\newcommand{\vEout}{\mathbf{E}_{\text{out}}}
\newcommand{\vF}{\mathbf{F}}
\newcommand{\vG}{\mathbf{G}}
\newcommand{\vS}{\mathbf{S}}
\newcommand{\vSig}{\mathbf{\Sigma}}
\newcommand{\vHatS}{\hat{\mathbf{S}}_\mathcal{D}}
\newcommand{\vHatSig}{{\mathbf{\Sigma}_{\hat{\mathbf{S}}, \mathcal{D}}}}
\newcommand{\vO}{\mathbf{0}}

\usepackage{bm} % Bold Greek letters

\author[1]{Tyler H. Chang}

\author[2]{Stefan M. Wild}

\affil[1]{Mathematics and Computer Science Division,
    Argonne National Laboratory,\\
    9700 S Cass Ave Bldg 240, Lemont, IL, USA 60439\\
    {\tt tchang@anl.gov}}

\affil[2]{Applied Mathematics and Computational Research Division,
    Lawrence Berkeley National Laboratory,\\
    1 Cyclotron Rd, Berkeley, CA, USA 94720\\
    {\tt wild@lbl.gov}}

\title{Designing a Framework for Solving Multiobjective Simulation Optimization
Problems}

\begin{document}

\maketitle

\begin{abstract}
Multiobjective simulation optimization (MOSO) problems are optimization
problems with multiple conflicting objectives, where evaluation of at least one
of the objectives depends on a black-box numerical code or real-world
experiment, which we refer to as a simulation.  
While an extensive body of research is dedicated to developing new algorithms
and methods for solving these and related problems, it is challenging and time
consuming to integrate these techniques into real world production-ready
solvers.  This is partly due to the diversity and complexity of modern
state-of-the-art MOSO algorithms and methods and partly due to the complexity
and specificity of many real-world problems and their corresponding computing
environments.  The complexity of this problem is only compounded when
introducing potentially complex and/or domain-specific surrogate modeling
techniques, problem formulations, design spaces, and data acquisition
functions.  
This paper carefully surveys the current state-of-the-art in MOSO
algorithms, techniques, and solvers; as well as problem types and computational
environments where MOSO is commonly applied.  We then present several key
challenges in the design of a Parallel Multiobjective Simulation Optimization
framework (ParMOO) and how they have been addressed.  Finally,
we provide two case studies demonstrating how customized ParMOO solvers can be
quickly built and deployed to solve real-world MOSO problems.
\end{abstract}


\section{Introduction and Motivation}

Multiobjective optimization problems (MOOPs) are optimization problems
involving two or more potentially conflicting objectives.  Such problems arise
in numerous fields of science, with examples such as multidisciplinary
engineering design \cite{ss2015,zhao2018}, scientific model calibration
\cite{bollapragada2020}, high-performance computing (HPC) library autotuning
\cite{chang2020c}, particle accelerator design \cite{neveu2023}, neural
network architecture search \cite{kandasamy2020,karl2022,parsa2020}, and
computational chemistry \cite{schweidtmann2018,shields2021}.  In such problems
the goal is to find a set of ``good'' points from a design space with respect
to a vector-valued objective function $\vF : \cX \rightarrow \R^o$.  Here,
$\cX$ is called the feasible design space and is typically assumed to be a
compact simply-bounded subregion of $\R^n$.  In this paper we assume each
objective $F_i(\vx)$ is bounded from below so that the feasible objective space
$\vF(\cX)$ is lower bounded in each dimension.  In the standard formulation,
one seeks to minimize each objective (i.e., component of $\vF$), a problem that
is written as
\begin{equation}
\label{eq:moop}
\min_{\vx\in\cX} \vF(\vx).
\end{equation}

Since it is typically not possible to find one $\vx \in \cX$ that
simultaneously minimizes all $o$ components of $\vF$, the solution to
\cref{eq:moop} is generally a set of design points and corresponding objective
values.  Rather than describe when an objective value corresponds to a solution
for \cref{eq:moop}, it is easier to begin by describing when an objective value
does {\it not} correspond to a solution.  An objective value
$\vF(\vx)\in\vF(\cX)$ is said to be dominated if there exists a $\vy \in \cX$
such that $\vF(\vy) \leq \vF(\vx)$.  Here we use the vector-inequality
``$\leq$'' to indicate that $\vF(\vy)$ is componentwise less than or equal to
$\vF(\vx)$ with strict inequality $F_i(\vy) < F_i(\vx)$ in at least one
component $i$.  Conversely, an objective value $\vF(\vx^\star)$ is said to be
nondominated in a set of objective values $\Omega\subset\R^o$ if for all
$\vF(\vy)\in\Omega$, $\vF(\vy) \not\leq \vF(\vx^\star)$.

If $\vF(\vx^\star)$ is nondominated in $\vF(\cX)$, then $\vF(\vx^\star)$ is
said to be Pareto optimal for \cref{eq:moop}, and $\vx^\star$ is said to be
efficient for \cref{eq:moop}.  The set of all Pareto optimal objective values
typically forms a $(o-1)$-dimensional trade-off surface, called the Pareto
front (along a subset of the boundary of $\vF(\cX)$).  For further details on
general MOOPs, we refer readers to the book by \cite{ehrgott2005}.  In
general, a design point and objective value pair $(\vx^\star, \vF(\vx^\star))$
are in the solution set for \cref{eq:moop} if they are efficient/Pareto
optimal.  In a typical multiobjective optimization application, a {\it decision
maker} uses a multiobjective optimization solver to identify a subset of all
Pareto optimal objective values and select one or more of these as the
preferred solution(s) based on domain expertise and personal preference.

A large class of real-world MOOPs consists of those derived from expensive
numerical simulations \cite{bollapragada2020,neveu2023,ss2015,zhao2018},
computer experiments \cite{chang2020c,kandasamy2020,karl2022,parsa2020}, and
real-world experiments \cite{elias2020,myers2016,schweidtmann2018,shields2021}.
We refer to such problems as multiobjective simulation optimization (MOSO)
problems.
Often, the simulations in a MOSO problem are black-box processes that are
expensive to evaluate and do not admit derivative information.  While the
objectives themselves may have components that are less expensive to evaluate
and/or admit derivative information, the MOSO problem as a whole is typically
treated as a black box and solved by using black-box optimization techniques
\cite{audet2017}.  Furthermore, due to the expected expense of evaluating the
simulations, 
the performance of MOSO algorithms
is typically measured in terms of solution accuracy
per number of total simulation evaluations.  In this paper we focus on
designing solvers for MOSO problems.  Consequently, here we use the convention
that ``simulation'' refers to an expensive black-box function that yields only
zeroth-order information.
As described later in this paper, solving MOSO problems does not always preclude the
usage of some partial derivative information.  However, we restrict
ourselves to techniques that account for the presence of a computationally
expensive black-box process somewhere in the computing chain, and we
work under the assumption that the complete gradient for all objectives and
constraints will not be available.

From an optimization researcher's perspective, solving MOSO problems
efficiently requires coordination between state-of-the-art and emerging
numerical techniques from a variety of fields of study, including
design-of-experiments, surrogate modeling, scalarization, uncertainty
quantification, and optimization algorithms. We elaborate on
many of these techniques in \Cref{sec:techniques}.

From a practitioner's perspective, solving real-world MOSO problems
involves more than numerical algorithms.  First, even within the narrow field
of MOSO problems, there can be significant variation in problem definitions.
This includes application and domain-specific modeling techniques, problem
structures, and feasible design spaces.  Second, it is important to consider
the resources and computing environment where the problem will be solved.  This
includes wet-lab environments and HPC resources, any of which may require
additional steps and technologies to efficiently operate.  We elaborate on
these challenges in \Cref{sec:challenges}.

Coordinating all of the above techniques from different fields of study into
usable solvers that can be deployed to solve real problems is a massive and
time-consuming undertaking, which often limits the rate at which novel methods
can be adopted.  The Parallel Multiobjective Simulation Optimization (ParMOO)
library implements a framework based on surrogate modeling for solving MOSO
problems, with the goal of providing an easy-to-use and modular interface for
implementing and deploying a wide variety of MOSO algorithms and emerging
techniques.  ParMOO was originally introduced in \cite{chang2023b}, with a
focus on the quality of the open-source software processes, intended user-base,
and usability.  This paper introduces the considerations that have continued to
drive ParMOO's design and development, and the MOSO algorithms and techniques
that ParMOO currently supports.  We begin by deconstructing surrogate-based
MOSO solvers into a system of simpler modules for performing subtasks that are
common across many solver architectures.  We then present an abstract
framework for integrating these modules, while reducing complexities.  To
demonstrate ParMOO's effectiveness in exploiting real-world problem structures,
we conclude by providing access to two computationally inexpensive benchmark
problems exhibiting some of the challenging features discussed above. These
problems are easy to test on, since the expensive simulations have been
replaced by machine learning models trained on real-world datasets.

The remainder of the paper is organized as follows.  In \Cref{sec:background}
we review existing techniques, problem types, trends, and software in MOSO, and
we describe the technical challenges and trade-offs involved in building a
general-purpose MOSO library that achieves all of our design goals.  In
\Cref{sec:framework} we introduce our framework for addressing these
challenges.
In \Cref{sec:dtlz2} we perform a brief performance study and validation on a
well-known multiobjective optimization test problem.
We then present two case studies that demonstrate the flexibility
and power of our framework: \Cref{sec:fayans} involves the calibration of a
physics model, and \Cref{sec:chem} involves the design of material in a wet-lab
environment.  We conclude in \Cref{sec:discussion} with a short discussion of
our results.

\section{Background and Design Principles}
\label{sec:background}

We now lay out the primary goals and requirements when designing a general
framework for implementing modern surrogate-based MOSO algorithms.  These goals
and requirements depend on an understanding of the current state-of-the-art
techniques and modern trends in MOSO, so we begin by reviewing these.

There are many different families of multiobjective optimization algorithms
based on where in the optimization process the decision maker expresses their
preference for a solution.  In this paper, we are interested in {\it a
posteriori} methods, where the decision maker expresses their preference after
solving the problem.  Other methods such as {\it a priori} methods and {\it
interactive} methods require different techniques.  They are summarized in the
survey by \cite{marler2004}.

It is useful to
further divide a posteriori MOSO methods into three broad categories that are common in
the literature: multiobjective evolutionary algorithms and
metaheuristics (MOEAs), which rely on random combinations and/or perturbations
to a so-called, ``pool,'' ``archive,'' or ``population'' of previously
evaluated design points; multiobjective direct search (MDS) methods, which rely
on (often deterministically) evaluating promising design points around one
or more current iterates in search of multiobjective descent directions; and
multiobjective Bayesian optimization (MBO), which focuses on minimizing a
multiobjective {\it acquisition function} to select design points for
evaluation.  In \Cref{sec:techniques}, we summarize into each of the
above methods, and the current representative algorithms for solving MOSO
problems.  Additionally, there are other methods -- mostly based on either
generic surrogate modeling or response surface methodology (RSM) -- which do
not fit cleanly into any of the above categories.  In \Cref{sec:techniques}, we
 explore these techniques as well.

\subsection{Fundamental Methods for Solving MOSO Problems}
\label{sec:techniques}

We now summarize the key techniques that are prevalent among
a posteriori MOSO algorithms, focusing on those that are most relevant to our
work; for a complete survey, see \cite{hunter2019}.

We begin by discussing scalarization, which is a classical technique that is
still prevalent in many modern MOSO algorithms and applications.  Scalarization
is the basis for almost all a priori methods, as it collapses a MOOP into a
single-objective subproblem that can be (approximately) solved with a
single-objective optimization solver.  In a priori methods, the decision maker
may provide some fixed set of scalarization functions, which are carefully
selected to target only the desired solution point(s), as pre-specified by the
decision maker.  However, in the context of a posteriori methods, a family of
multiple scalarization functions can be used to produce a discrete
approximation to the efficient set and/or Pareto front.

A scalarization is defined by a scalarization function $Q : \R^o \rightarrow
\R$, such that
\begin{equation*}
\arg \min_{\vx\in\cX} Q(\vF(\vx)) \subset \arg \min_{\vx\in\cX} \vF(\vx).
\end{equation*}
The scalarized problem $\min_{\vx\in\cX} Q(\vF(\vx))$ can be solved by using
any derivative-free optimization method.
The most commonly used scalarization technique is the weighted sum
scalarization \citep[Ch.~3]{ehrgott2005}, where the objective is collapsed
via a weighted average.
Other common scalarization techniques include the epsilon-constraint method
\citep[Ch.~4]{ehrgott2005}, Pascoletti-Serafini scalarization
\cite{eichfelder2009}, the reference point method \cite{wierzbicki1999},
quadratic scalarization schemes \cite{dandurand2016},
and (augmented) weighted Chebyshev scalarization \cite{steuer1983}.

In order to provide an approximation to the complete Pareto front,
each of the above-mentioned approaches would need to be combined with some {\it
adaptive scheme} for sweeping through different members in a family of
scalarization functions, thus producing numerous solution points covering the
Pareto front \cite{das1998,deshpande2016,eichfelder2009,laumanns2006}.
Although it is a classical technique, scalarization is still widely
used in MOEAs \cite{deb2013,jain2013}, MDS \cite{audet2010,deshpande2016},
and MBO methods \cite{knowles2006}.  Even algorithms that do not explicitly
rely upon scalarization still use some of these scalarizations indirectly.  For
example, the MDS method MOIF \cite{cocchi2018} and the line search technique
DFMO \cite{liuzzi2016} both select a multiobjective steepest descent direction
that could be viewed through the lens of scalarization.

Many newer MOSO algorithms (particularly in MOEAs) have replaced the
scalarization function with an {\it indicator function}, which focuses on
assessing the quality of the Pareto front approximation when adding a new point
$\vF(\vx)$.  Given a dataset $\cD = \{(\vx_i, \vF(\vx_i))\}_{i=1}^N$ of
previously evaluated design points and objective values, an {\it indicator
function} $I(\cD)$ calculates the quality of the Pareto front approximation
given by the nondominated objective values in $\cD$.  Designing a good indicator function is also nontrivial, but it is
generally agreed that a good indicator should produce solution points that are
on or close to the true Pareto front, offer good coverage of the entire Pareto
front, and are well distributed \cite{audet2020,chang2022a}.  By choosing a
scalarization function $Q_\cD(\vx) = I(\cD \cup \{(\vx, \vF(\vx))\}$ or
$Q_\cD(\vx) = -I(\cD \cup \{(\vx, \vF(\vx))\}$, the indicator can play the part
of both the scalarization function and its adaptive scheme.  The key difference
here is that an indicator-function-driven scalarization depends upon the
dataset $\cD$ in addition to the value of $\vF(\vx)$.  One of the most commonly
used indicators in MOSO algorithms is the hypervolume improvement indicator,
which measures the amount of hypervolume between a potential objective value
and previously observed objective values \cite{bringmann2013,shang2020}.  This
indicator is used as the basis for many multiobjective algorithms, such as
Archived Multiobjective Simulated Annealing (AMOSA) \cite{bandyopadhyay2008}.

Since each scalarization function reduces a MOSO problem to a single-objective
subproblem, this technique must be coupled with a single-objective solver.
Such a solver could implement a single-objective evolutionary algorithm,
heuristic, direct search method, Bayesian optimization solver, or any other
single-objective approach.  It is also possible to apply a different
MOSO solver (or a single iteration thereof) to the MOSO problem, then select
the solutions that best minimize the scalarization a posteriori.  For a recent
survey of derivative-free optimization methods that could be used in the above
context, see \cite{larson2019}.

For MOSO problems where simulation evaluations are expensive, the total problem
expense can get out of hand if too many simulation evaluations are required.
Therefore, many modern MOSO solvers reduce the need for true simulation
evaluations by combining one or more of the above-mentioned multiobjective
algorithms with a computationally cheaper surrogate model
\cite{bradford2018,chang2022a,muller2017}.  Here, each component of an
expensive-to-evaluate objective $\vF(\vx) = (F_1(\vx),\ldots,F_o(\vx))$ is
modeled by a computationally cheaper surrogate function
${\hat \vF}_\cD = \left({\hat F}_1, \ldots, {\hat F}_o\right)$, based on the
current dataset $\cD$, and such that
${\hat F}_1 \approx F_1,\ldots,{\hat F}_o\approx F_o$.  Several ``candidate
solutions'' can then be suggested by finding points that are efficient/Pareto
optimal for the surrogate problem
\begin{equation}
\min_{\vx\in\cX} \big({\hat F}_1(\vx), \ldots, {\hat F}_o(\vx)\big).
\label{eq:surrogate}
\end{equation}
When $\vF$ is sufficiently more expensive to evaluate than ${\hat \vF}_\cD$ in
terms of computational time and resources, the cost of finding (approximately)
efficient/Pareto optimal points for \cref{eq:surrogate} can be much less than
the cost of finding such points for \cref{eq:moop}.  In practice, these
surrogates could be any of a variety of classical approximation or machine
learning \cite{bouhlel2019} models.

The main challenge when using surrogate models is that the surrogate models
must be asymptotically accurate for the optimization algorithm to convergence
to an accurate solution.  Typically, the accuracy of the surrogate depends on
the geometry and density of the dataset $\cD$ in $\cX$, the componentwise
smoothness of the underlying function $\vF$, and the choice of surrogate model
\cite{cheney2009}.  In most applications, the smoothness of $\vF$ is an
unknown but fixed quantity.  Therefore, we must rely on either the density and
geometry of $\cD$ or {\it uncertainty information} returned from the surrogate
models in order to ensure acceptable surrogate accuracy.

In the context of model-based derivative-free optimization (DFO), one focuses
on the local properties of $\cD$ around the current iterate $\vx^{(k)}$.  As
more design points are evaluated during the course of the optimization
algorithm, the density of $\cD$ rapidly increases near local minima.
Therefore, the focus is on maintaining an acceptable conditioning for the
surrogate modeling problem in the neighborhood of such minima, which is
determined by the local geometry of the design points in $\cD$
\cite{conn2008}.  This geometry can be maintained by occasional requiring the
optimization algorithm to evaluate model-improving design points instead of
strictly model minimizing points \cite{conn2009}.  Since these techniques,
focus on local quality of the surrogate model, the surrogate problem
\Cref{eq:surrogate} must be constrained with a local trust region (LTR) so that
only candidate solutions in the LTR are suggested by the optimization
algorithm.  There are currently a few surrogate-driven MOSO algorithms that
exactly follow this approach \cite{thomann2019,berkemeier2021}, but none with
software to our knowledge.  Other existing algorithms use aspects of this
approach or are similar in spirit \cite{chang2022a,ryu2014}.  From the
single-objective setting, this ``model-based DFO'' approach is known to produce
fast converging highly scalable algorithms, but only guarantees local
convergence \cite{ragonneau2021}.  However, this approach can be coupled with
global search or multi-start techniques if global accuracy is desired
\cite{larson2018}.

On the other hand, Bayesian optimization leverages uncertainty models specific
to surrogates.  Certain surrogates, such as Gaussian processes, admit an
uncertainty function $\vSig_{{\hat \vF}, \cD} : \cX \rightarrow \R_+^o$ such
that $\vSig_{{\hat \vF}, \cD}(\vx) \propto \left| \vF(\vx) - {\hat \vF_\cD}(\vx)
\right|$, (where `$\propto$' can be understood as meaning ``componentwise
proportional'') \cite{garnett2023}.   In MBO, the scalarization/indicator
functions are augmented with the value of $\vSig_{{\hat \vF}, \cD}$ to assess
the overall utility of evaluating a potential candidate design point $\vx$, in
terms of both improving the Pareto front approximation and maintaining global
accuracy of the surrogate model \cite{emmerich2016,feliot2016}.  Notably, the
ParEGO algorithm uses the uncertainty information from a Gaussian process model
to calculate the expected improvement in the augmented Chebyshev scalarization
\cite{knowles2006}, and the (q)EHVI algorithm uses (a quasi-stochastic
approximation to) the expected improvement in the hypervolume indicator
\cite{daulton2020}.  MBO produces globally convergent methods, but with the
drawback of a slower rate of convergence, particularly in high-dimensional
design spaces.  Therefore, many methods aimed at high-dimensional design spaces
use localization strategies such as LTRs, similarly to the model-based DFO
methods above \cite{eriksson2019}.

An alternative framework from the statistics literature is multiobjective RSM.
RSM begins with a ``search'' step, which generates and evaluates a large
design-of-experiments \cite{garud2017} or quasi-random sample \cite{roy2023}
to explore the design space and collect an initial dataset for surrogate modeling.
Once enough data has been collected to ensure sufficiently accurate surrogate modeling, the
typical multiobjective RSM approach applies a polynomial surrogate model,
multiple scalarizations, and an optimization procedure to solve all scalarized
surrogate problems \citep[Ch.~7]{myers2016}.  Several existing multiobjective
surrogate-modeling algorithms use variations or multiple iterations of this
approach, including VTMOP \cite{chang2022a}, SOCEMO \cite{muller2017}, and
PAWS \cite{ryu2014}.  Multiobjective RSM can be efficient in many-objective
settings since a large portion of RSM's computational cost comes from the
initial search step, whose results can be shared across all scalarizations.
Therefore, the incremental cost of solving for many objectives by applying many
scalarizations is not much greater than the cost of solving for a single
scalarization by RSM.  The drawback of this approach is that in
high-dimensional spaces, it may not be feasible to generate a sufficiently
large design-of-experiments to guarantee global surrogate accuracy.  Still, all
of the above methods also require an initial surrogate-modeling dataset, which
is often gathered by similar techniques as in RSM; see, e.g., the search
techniques used in the Bayesian optimization library BoTorch \cite{balandat2020}.

In practice, the lines between the types of algorithms discussed above are
often blurred.  Furthermore, not all MOSO algorithms use all of the methods
described above.  In contrast, ParMOO focuses on techniques that utilize:
\begin{itemize}
\item An acquisition function or similar technique for setting targets in the
objective space;
\item An optimization solver (or single iteration thereof) for generating the
next iterate based on the acquisition function;
\item A surrogate model for approximating expensive functions and either an
uncertainty function or model improvement procedure for maintaining its
accuracy; and
\item A search technique for exploring the design space and generating the
initial samples for surrogate modeling.
\end{itemize}

Prior to the release of ParMOO, there was no open-source
software for implementing generic MOSO solvers that adhere to the above
structure.  \cite{chang2023b} focused on the usage of ParMOO
rather than the underlying algorithm, framework, and techniques supported.  In this paper, we present the design strategies and
architecture used for solving MOSO problems via any algorithm that falls into
the above framework, with a focus on how such combinations can be most effectively applied
in domain-specific applications.

\subsubsection{Other Techniques in MOSO.}
\label{sec:other-techniques}

As previously stated, not all MOSO algorithms utilize all of the techniques
discussed in \Cref{sec:techniques}.

One of the biggest departures from \Cref{sec:techniques} is that not all MOSO
algorithms use a technique that can be cast under the umbrella of an
acquisition function.  However, most require a technique for selecting the
point or set of points to iterate from.  For example, the MDS method MODIR
\cite{campana2018} generalizes the single-objective DIRECT algorithm by way of
a multiobjective identification procedure.  The MDS methods
DMulti-MADS \cite{bigeon2020}, DMS \cite{custodio2011}, and MultiGLODS
\cite{custodio2017} use some combination of indicator functions and/or
trust-region radii to select each iterate.  In the context of MOEAs, a
dominance-based sorting metric is typically applied during population selection
\citep[Ch.~2.3.3]{abraham2005}.  Such an approach is taken in the well-known
MOEA NSGA-II \cite{deb2002a}.  While many of these selection metrics could be
calculated from the information available to an acquisition function (in
particular, the dataset $\cD$), these methods represent a significant departure
from the techniques discussed in \Cref{sec:techniques}.
Additionally, it is worth noting that many of the algorithms listed above do
not rely on surrogate modeling as a core feature, which is another significant
departure from the techniques covered toward the end of \Cref{sec:techniques}.

While not all encompassing, the class of techniques covered in
\Cref{sec:techniques} give the core ingredients for many
MOSO algorithms.  Algorithms of this nature are the primary focus of this
paper.

\subsection{Design Challenges when Building MOSO Solvers}
\label{sec:challenges}

We now explore the key design challenges addressed by ParMOO.
ParMOO is designed to use the techniques laid out in \Cref{sec:techniques},
specifically, acquisition/indicator/scalarization functions, surrogate
modeling, and design-of-experiments, to solve a diverse set of scientific MOSO
applications.  Therefore, we focus on challenges related to these
techniques and how they can be deployed in real-world applications and
scientific computing environments.

\subsubsection{State-of-the-Art and Domain-Specific Techniques.}

First, as highlighted in \Cref{sec:techniques}, state-of-the-art
surrogate-driven MOSO solvers require integrating several different techniques,
including acquisition functions, optimization solvers, surrogate models, and
design-of-experiments.  This challenge is compounded when we consider the
inclusion of domain- and problem-specific techniques such as machine learning
surrogates with physics-informed constraints \cite{karniadakis2021} or
acquisition functions designed to target extremely rare events
\cite{pickering2022}.

In ParMOO, we seek to support most such techniques, including domain-specific
methods that we may not yet be aware of.  This requires a clean application
programmer's interface (API) for users to implement new, domain-specific
methods.  This, in turn, requires a standardization of APIs for each individual
method and clear definition of how these methods will interact, which may not
be easy to define for all of the techniques in \Cref{sec:techniques}.

\subsubsection{Structured Problem Formulations.}

Beyond using domain- and problem-specific algorithms, one can go a step farther
and explicitly model domain- and problem-specific structures in the MOSO
problem formulation.  This includes modeling composite problem definitions, as
described in \cite{astudillo2021,khan2018,larson2023}, where the objective can be
expressed as the composition of one or more functions.  In particular, we can
modify the MOOP problem definition from \Cref{eq:moop}, using
\begin{equation}
\min_{\vx \in \cX} \vF(\vx, \vS(\vx))
\label{eq:moso}
\end{equation}
where $\vF : \cX \times \R^m \rightarrow \R^o$ and $\vS : \cX \rightarrow
\R^m$.  In this formulation, it is assumed that $\vF$ is an {\it algebraic}
function that is cheap to evaluate and whose gradient function may be known,
while $\vS$ is a computationally expensive simulation function, as previously
discussed.  This formulation decouples the ``simulation'' aspect of the
computation from the objective function, allowing for the usage of structured
non black-box optimization solvers.

One common example for $\vF$ would be a sum-of-squares function
$F_i(\vx, \vS(\vx)) = \sum_{j\in J} S_j(\vx)^2$, for one or more component
functions $F_i$ and an index set $J$.  In the single-objective case, this
structure can be exploited by Gauss-Newton inspired methods to achieve
super-linear convergence rates \cite{zhang2012}, as is the case in the
derivative-free least-squares solver POUNDERS \cite{wild2017}.  While the
theoretical advantage may not be as strong, explicitly modeling the composition
in \Cref{eq:moso} typically leads to convergence with fewer simulation
evaluations for a variety of other algebraic $\vF$-functions
\cite{astudillo2021,khan2018,larson2023}.

Another structure of $\vF$ that is specific to the multiobjective setting
is a heterogeneous problem definition \cite{thomann2019}, where one or more
$F_i$ is an algebraic function with no dependence on an expensive simulation
output, but other $F_i$ are either the direct outputs or calculated from a
simulation output.  This sort of problem definition has been acknowledged in
previous works \cite{chang2022a,thomann2019}, but we are not aware of a formal
analysis.  The typical approach in the context of surrogate-based MOSO, is to
apply surrogate models to the components that depend upon simulation outputs
and use the algebraic formulations instead whenever they are available.  This
approach lends well to the abstraction in \Cref{eq:moso}.

In ParMOO, we support problem definitions of the form \Cref{eq:moso}, which
is a slight departure from the standard black-box optimization methodology.
However, it is important to note that this formulation still involves a
black-box process, and therefore cannot be solved with purely gradient-based
methods.

\subsubsection{Changes to the Design Space -- Mixed Variables and Constraints.}
\label{sec:embeddings-constraints}

Thus far, we have assumed that $\cX$ is a compact, simply-bounded subregion of
$\R^n$.  However, widely used MOSO software often support broad categories
of design variables and constraints, including various combinations
of real (continuous-valued), discrete (integer-valued), and/or categorical
(nonordinal) variables, as well as linear and nonlinear constraints
\cite{balandat2020,benitezhidalgo2019,biscani2020,blank2020,kandasamy2020,
ledigabel2011}.

For mixed-type design variables, it is important to
note that changing the domain $\cX$ can result in a completely different
classification of problem requiring different techniques.  For example, when
$\cX$ consists of exclusively binary decision variables, a multiobjective
ranking and selection algorithm would generally be used \cite{feldman2018},
which significantly differs from the continuous MOSO algorithms described in
\Cref{sec:techniques}.  However, many MOSO problems that are largely continuous
have a small number of discrete variables, which is something that can be
reasonably supported without significantly changing the kind of techniques
used.

Since we are focused on methods that were designed for continuous MOSO, the
primary challenge is in handling the discrete design variables.  For MOEAs,
this would be supported by providing a mixed-variable mating procedure to any
existing MOEA \cite{blank2020}.  For MBO and many MDS techniques it could be
sufficient to define a distance, kernel, or neighborhood
function that supports integer and categorical variables \cite{audet2023}.
For a generic surrogate-based MOSO algorithm, however, one needs a function $E
: x_{i'} \rightarrow [0,1]^\ell$ for embedding the discrete design variable
$x_{i'}$ into a continuous latent space $[0,1]^\ell$, where $\ell$ is the
dimension of the embedding.  Perhaps the most common example of a latent space
embedding is a one-hot encoding, which maps each category to either a $0$ or
$1$, and assigns a value $1$ if and only if the corresponding $x_{i'}$
represents the corresponding category.  This latent space is then modeled by
the surrogate, which can be used to relax the embedded design variables to
a continuous representation.  This approach has been widely used in existing surrogate modeling
software \cite{bouhlel2019,saves2023}.  For any of the above techniques, it
can be difficult to scale generic techniques for problems with a large number
of categorical design variables.  Therefore, many domain-specific software
packages and algorithms utilize domain-specific design space embedder or
distance functions \cite{kandasamy2020,shields2021}.  As an example in the
context of material design, the Bayesian optimization solver EDBO
\cite{shields2021} uses the molecular descriptor calculator Mordred
\cite{moriwaki2018} to calculate a three-dimensional latent space encoding
molecule networks, and then solves a Bayesian optimization problem in this latent
space.

Similarly, the introduction of linear and nonlinear constraints requires us to
relax the previous assumption that $\cX$ is simply bounded.  In practice, such
constraints could be many combinations of known (e.g., part of the problem
formulation) or hidden (e.g., discovered at runtime), algebraic or dependent
upon simulation outputs, relaxable in the sense that the simulation can still
be evaluated when the constraint is violated, and quantifiable
in the sense that we can quantify how much it has been violated or
not \cite{ledigabel2015}.  Clearly, we will not be able to support
all combinations of the above, however, we seek to support as many
combinations as reasonably possible.

Given the composite objective function definition in \Cref{eq:moso}, it is
natural and simplifying to mirror this with a composite constraint definition
$\vG(\vx, \vS(\vx)) \leqq \vO$, where ``$\leqq$'' denotes componentwise less than
or equal to.  In this sense, we will support both simulation-based and
algebraic inequality constraints and immediately preclude hidden constraints.
Then, we can decompose the feasible region $\cX$ into a simply-bounded region
$\cB$ (with $n$ dimensions and possibly mixed variable types), where the bound
constraints are unrelaxable, and the relaxable nonlinear constraints defined by
$\vG$.  This produces ParMOO's final problem formulation
\begin{equation}
\min_{\vx \in \cB} \vF(\vx, \vS(\vx)) \text{ subject to } \vG(\vx, \vS(\vx))
\leqq \vO.
\label{eq:parmoo}
\end{equation}

For generic nonlinear black-box constraints, there have been several approaches
in the single-objective optimization literature for handling the nonlinear
constraints $\vG$, including barrier functions \cite{audet2009} and augmented
Lagrangian approaches \cite{lewis2002}.  Both of these approaches can be
thought of as producing a penalty function describing the constraint violation,
which can be added to the objective to encourage feasibility.  However, the
methods for achieving this differ.  To extend these methods to the
multiobjective case, it is sufficient to apply a cumulative penalty function to
all objectives \cite{cocchi2020}.  Considering the options, barrier functions
tend be more flexible to apply without significant modification to the
optimization algorithm.  Among the barrier methods, the two common methods are
an extreme barrier, which adds an infinite penalty any time a constraint is
violated, and a progressive barrier, which adds a penalty based on the
constraint violation and progressively increases this penalty over time.

In ParMOO, we seek to support a limited number of mixed variables and
relaxable nonlinear constraints $\vG$.  We pursue this through embeddings
and barrier functions, which requires limited modification to the
techniques outlined in \Cref{sec:techniques}.

\subsubsection{Parallel and Other Computing Environments.}

Another challenge in MOSO is achieving efficient resource utilization in
diverse computing environments.  Many existing optimization solvers provide
parallelism through a single paradigm (e.g., Python multiprocessing, OpenMP, or
MPI) \cite{blank2020,chang2022a,ledigabel2011}.
However, this may not be sufficient when dealing with diverse computing
environments including distributed systems, HPCs, laboratory clusters, and even
wet-lab environments, any one of which may not support one or more common
parallel computing paradigms.

In order to address this complexity, many platforms \cite{chard2020,elias2020}
and libraries/frameworks \cite{hudson2022a,hudson2022b} have been created for
coordinating parallel evaluations, complex simulations, data storage, and
heterogeneous resource utilization.  Additionally, domain-specific tools for
managing simulation, modeling, and optimization interplay also manage
distribution over parallel resources \cite{kolonay2011}.  However, effectively
utilizing these platforms is nontrivial if the optimization software was not
originally designed with mindfulness for these paradigms
\cite{chang2020b,raghunath2017}.

In ParMOO, we have designed with parallelism and simulation environment
flexibility in mind, being sure to decouple simulation evaluation from the
solvers and techniques.

\subsubsection{Maintainability and Usability.}

A recent movement in the general scientific computing space is to reduce
technical debt and improve scientist productivity through better scientific
software development practices \cite{heroux2020}.
Many scientific software libraries are adopting such practices
\cite{hudson2022a,ragonneau2021}, and several existing multiobjective
optimization libraries have recently been refactored to improve maintainability
\cite{audet2022,benitezhidalgo2019}.

One of our key goals in the design of ParMOO has been to build a software
package that provides an API for implementing the techniques discussed in
the previous two sections, without sacrificing maintainability and usability.

\subsubsection{Design Goals.}
\label{sec:des-goals}

In summary to the challenges and requirements laid out in this section, we now
list five key goals in the design of ParMOO.
\begin{enumerate}
\item[1.]
We want to provide a highly customizable software framework for building and
deploying surrogate-based MOSO solvers.
\item[2.]
We want to give the optimization solver access to any available structure
in how the simulation outputs are used to define the objectives.
\item[3.]
We want to be flexible in our support for a wide variety of problem types,
including mixed variables and constraints.
\item[4.]
We want to make our framework easy to deploy in a wide variety of scientific
workflows.
\item[5.]
Our software framework and workflow must be easy to use, maintain, and
extend.
\end{enumerate}

Goals 1, 2, and 3 require careful thought in terms of what kinds of MOSO
algorithms, problem types, and structures are supported and how these
techniques will interact.  Goals 4 and 5 are engineering problems, which
constrain the breadth of methods that we can reasonably support when addressing
Goals 1, 2, and 3.  Although it constrains our
decision making on the techniques we can use, ParMOO's workflow and
general project structure (some of the main components of Goal 5) have
previously been extensively discussed in \cite{chang2023b}.  Therefore, we limit our discussion regarding Goal 5 to
where it affects our decision making for the other design goals.

\subsection{Existing MOSO Software}
\label{sec:sw}

In addition to the software laid out below,
we acknowledge the existence of other special-purpose solvers for
domain-specific variations of the MOSO problem.  In this work, however, we do
not consider solvers that are tied to a specific application, such as neural
architecture search \cite{parsa2020} and chemical experiment design
\cite{shields2021}, or solvers that target other variations of the problem
\cref{eq:moop}, including online multiobjective optimization \cite{mannor2014}
and multiobjective reinforcement learning \cite{hayes2022}.  We also
acknowledge that a similarity exists between surrogate-based MOSO and active
learning \cite{sapsis2022}, which we consider to be a more generic problem
than surrogate-based MOSO.

In this section, we exclusively list open-source, production-quality software
packages for solving MOSO problems.  The term ``production quality'' here focuses on software implementations that are advanced
enough for a non-expert to use in a non-academic application.  Hallmarks of
such software include some combination of detailed user or API documentation,
publication in a software journal, and a large user base.  Although often not distinguished, it is also
important to clarify that not every well-known MOSO algorithm has a one-to-one
mapping with such a software implementation.

Much of the available software
implements MOEAs.  Widely used MOEA-based
solvers include ParEGO \cite{knowles2006}, which also integrates Gaussian
process surrogate modeling, and SPEA2 \cite{zitzler2001}.  There are also many
libraries, including Platypus \cite{hadka2015}, pymoo \cite{blank2020},
PlatEMO \cite{tian2017}, jMetal/jMetalPy
\cite{durillo2011,benitezhidalgo2019}, and pagmo/pygmo \cite{biscani2020}.
The majority of these libraries implement all of the well-known MOEAs, such as
NSGA-II \cite{deb2002a} and NSGA-III \cite{deb2013,jain2013}.
Although not specific to MOOPs,
the framework DEAP \cite{fortain2012} is also frequently used for implementing
distributed evolutionary algorithms in Python, including distributed MOEAs.

Packages with MDS methods include the following.  The DFO-lib
\cite{liuzzi2024} contains serial Fortran implementations of MODIR and DFMO and
both Matlab and Python implementations of MOIF.  Serial and parallel Fortran
solvers are distributed in VTMOP \cite{chang2022a}, which utilizes RSM together
with adaptive weightings and trust-region methods.  NOMAD is a
parallel-capable C++ library of industrial-grade single- and multiobjective
direct search methods, which until recently only contained the biobjective
solver BiMADS \cite{ledigabel2011}.  In a newer release of NOMAD v4
\cite{audet2022}, an implementation of DMulti-MADS was added as well.  BoostDFO
\cite{tavares2021} is a MATLAB library containing parallel implementations of
single- and multiobjective direct search solvers, including DMS and MultiGLODS.
PyMOSO \cite{cooper2020} is a Python framework that is targeted primarily at
integer-valued problems.

In recent years several Bayesian optimization libraries and frameworks have
emerged, supporting MBO.  The Python library Dragonfly \cite{kandasamy2020},
was designed for solving neural architecture search problems but can be applied
to generic MBO problems.  Perhaps most relevant to our goals in this paper is
the BoTorch \cite{balandat2020} framework for implementing and deploying
parallel, production-ready, generic Bayesian optimization (including
multiobjective) algorithms while layering over the automatic differentiation
PyTorch framework \cite{paszke2019}.

To summarize this section, we present a list of the solvers and libraries
discussed thus far in \Cref{tab:moo-sw}.  This table
represents the relevance of these libraries in relationship to ParMOO's design
goals.  Therefore, we focus on the features offered by these libraries that
relate to Design Goals 1--4 from \Cref{sec:des-goals}.  Notably, support for
domain- and application-specific problem definitions (related to Goal 2) is not
shown since, of the methods listed, only BoTorch and ParMOO support this kind of
problem formulation at this time.
We acknowledge that many of these software packages are under active
development and may add features in the future, so we base our classifications
by the documented features at the time of this publication.

\begin{table}[h]
\begin{center}
\begin{tabular}{c|cccccc}
Name   & Type & Language & Method & Constr. & Var.\ Types & Surrogates\\
\hline \hline
NOMAD v4  & L & C++     & MDS     & yes  & mixed & yes  \\
BoostDFO  & L & Matlab  & MDS     & some & real  & yes  \\
BoTorch   & L & Python  & MBO     & yes  & mixed & yes  \\
DEAP      & L & Python  & MOEA    & yes  & mixed & no   \\
DESDEO    & L & Python  & any     & yes  & real  & yes  \\
Dragonfly & L & Python  & MBO     & yes  & mixed & yes  \\
jMetal/jMetalPy    & L & Java/Python    & MOEA    & yes  & mixed & no   \\
DFO-lib   & L & Fortran/Python/Matlab & MDS     & some   & real or int  & no   \\
ParEGO    & S & C       & MOEA/MBO& no   & real  & yes  \\
pagmo     & L & C++     & MOEA    & some & mixed & no   \\
ParMOO    & L & Python  & MDS/MBO & yes  & mixed & yes  \\
PlatEMO   & L & Matlab  & MOEA    & some & mixed & some \\
Platypus  & L & Python  & MOEA    & yes  & mixed & no   \\
pygmo     & L & Python  & MOEA    & some & mixed & no   \\
pymoo     & L & Python  & MOEA    & some & mixed & no   \\
PyMOSO    & L & Python  & MDS     & yes  & int   & no   \\
SPEA2     & S & C       & MOEA    & no   & real  & no   \\
VTMOP     & S & Fortran & MDS     & no   & real  & yes  \\
\end{tabular}
\vskip 12pt
\end{center}
\caption{
General-purpose open-source MOSO software.  For each software package,
the columns are labeled as follows:
``Type'' indicates whether this is an individual solver (S) or
library of solvers (L);
``Language'' specifies the primary development language;
``Method'' classifies each software as primarily using MOEA, MDS, or MBO;
``Constr.'' indicates whether nonlinear constraints are supported;
``Var.\ Types'' indicates the types of variables supported;
and ``Surrogates'' indicates whether surrogate modeling is used.
}
\label{tab:moo-sw}
\end{table}

\section{A Framework for Multiobjective Simulation Optimization}
\label{sec:framework}
We now outline our framework for solving MOSO problems.  In particular, we
describe how this framework addresses the design goals and challenges laid out
in \Cref{sec:des-goals}.
The framework is implemented in the software library ParMOO \cite{chang2023b},
which can be obtained by following the installation instructions in its online
documentation \cite{chang2023a}.

\subsection{An Object-Oriented Modular Design}
\label{sec:oop}

Despite the level of complexity necessary to achieve Goal 1, the biggest
challenge of this work is maintaining a usable, maintainable, and extensible
API while managing the complexity of these interacting techniques (Goal
5).  To do so, we focus on a modular design, where each module addresses a
particular component of the larger problem through a common interface.
The natural paradigm for implementing such a framework is object-oriented
programming (OOP).  To this end, we have implemented abstract base classes
(ABCs) for each of the solver components.
Referring back to \Cref{sec:techniques}, the four primary ABCs for ParMOO are
\begin{itemize}
\item the {\tt GlobalSearch} class, which abstracts the API for utilizing
design-of-experiments and other sampling techniques;
\item the {\tt SurrogateFunction} class, which abstracts the API for
fitting, updating, and evaluating a generic surrogate function, plus either
evaluating the uncertainty function or suggesting model improvement points, as
applicable;
\item the {\tt AcquisitionFunction} class, which abstracts the API for
implementing
an acquisition function (including scalarization and indicator functions) and
selecting the next iterate (if applicable); and
\item the {\tt SurrogateOptimizer} class, which accepts a scalarized,
differentiable or non-differentiable function and uses an appropriate method to
produce the next iterate or (if applicable) call for a model improvement.
\end{itemize}
To allow users to utilize ParMOO without implementing each of these
techniques, we provide a built-in library of implementations for several
standard techniques from each of the categories.

At the center of all of this, the {\tt MOOP} class stores a potentially complex
MOSO problem definition of the form \Cref{eq:parmoo} and coordinates
implementations of the ABCs to solve a MOSO problem.  This is achieved using
the ``builder'' design pattern \cite{gamma1995}.
The {\tt MOOP} class itself is also highly modular and extensible, with clearly
documented methods for performing
both private and public tasks.
Developers can easily overwrite certain functionalities, such as the method
that distributes expensive simulation evaluations, in order to adapt to novel
architectures and use cases.  In this way, the modular OOP framework also
addresses goal 5.

The downside of this approach is that the object-oriented API adds an
additional level of complexity, similar to that of learning a new modeling
language, such as JuMP \cite{dunning2017} or Pyomo \cite{hart2017}.  It also
requires users to become aware of the basic components of a surrogate-based
multiobjective solver, in order to correctly utilize our library of components
and framework.  Ultimately, this mental overhead is significantly greater than
that required to learn to properly utilize other Python libraries of
multiobjective solvers, such as pymoo \cite{blank2020} and pagmo/pygmo
\cite{biscani2020}.  However, we believe that this approach is a necessary
compromise in order to achieve Goals 1--3 as laid out in \Cref{sec:des-goals}.

In this way, ParMOO is not comparable to many of the libraries and individual
solvers listed in \Cref{tab:moo-sw}.  Of these software packages, ParMOO is
perhaps most comparable to DEAP \cite{fortain2012} and BoTorch
\cite{balandat2020}, which implement similar frameworks for solving slightly
different classes of problems using different classes of techniques.

\subsection{Distinguishing between Simulations, Objectives, and Constraints}
\label{sec:composite}

\begin{figure}[ht]
\begin{center}
\includegraphics[width=0.28\textwidth]{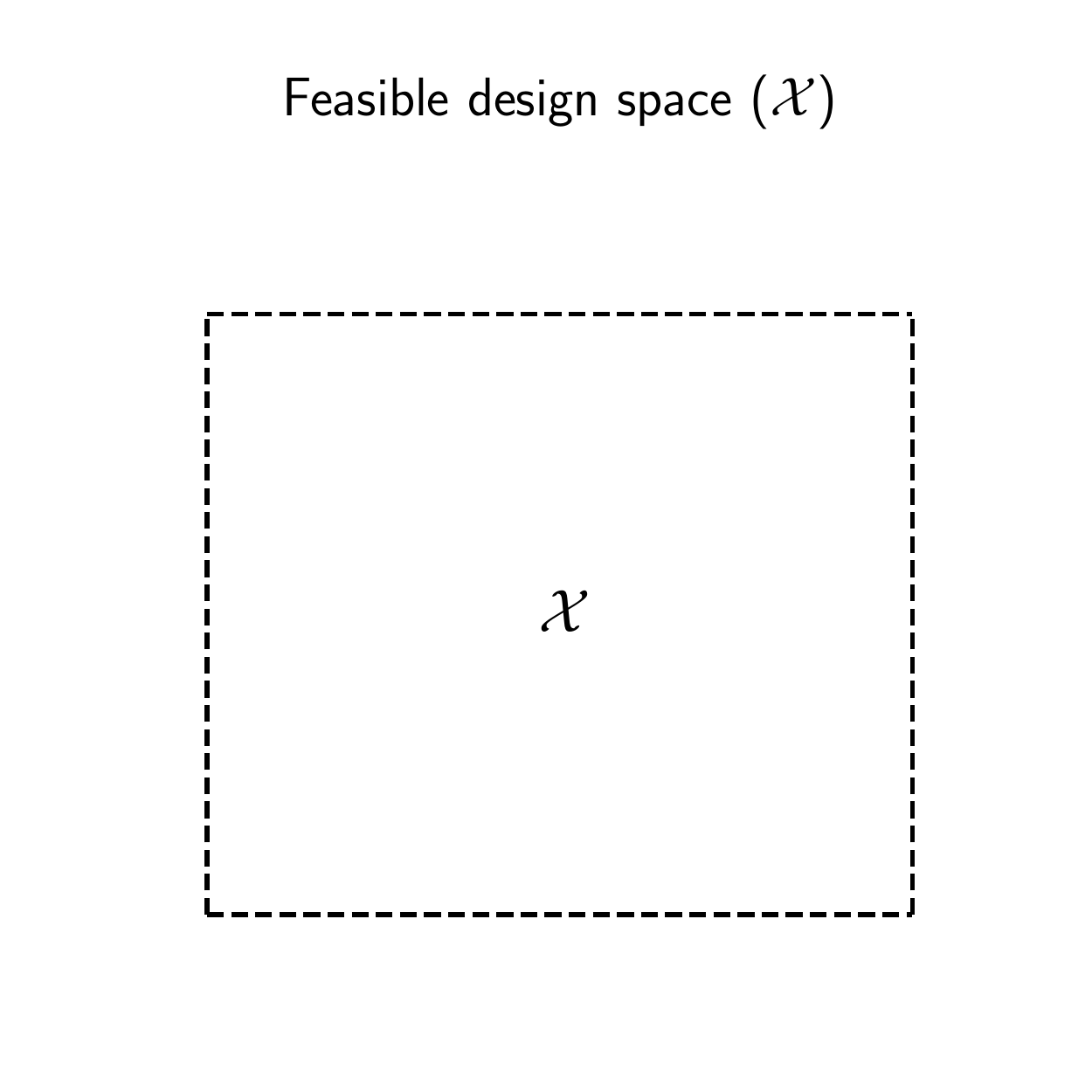}
$\begin{array}{c} \vS \\ \longrightarrow \\
\\ \\ \\ \\ \\ \\
\end{array}$
\includegraphics[width=0.22\textwidth]{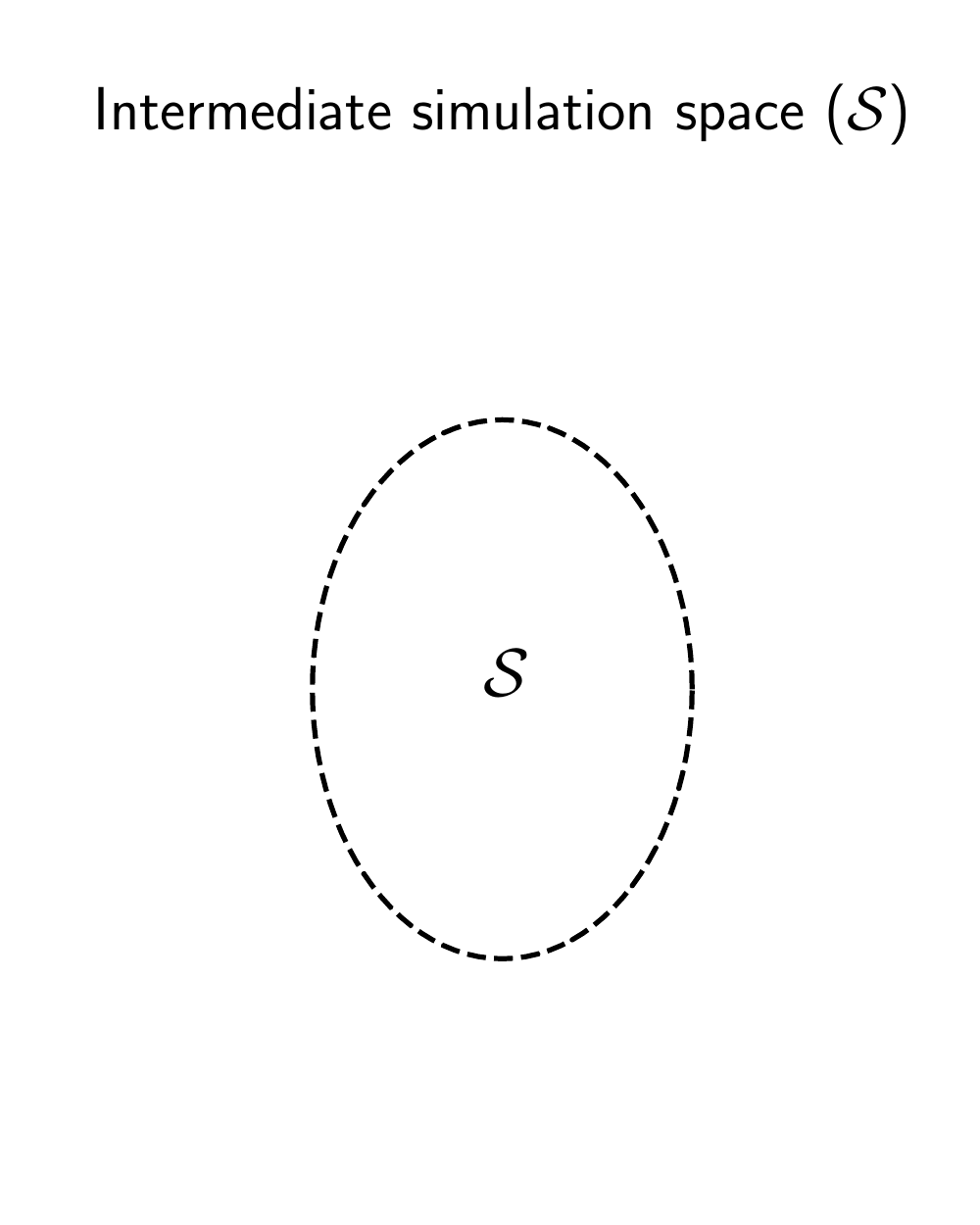}
$\begin{array}{c} \vF \\ \longrightarrow \\
\\ \\ \\ \\ \\ \\
\end{array}$
\includegraphics[width=0.28\textwidth]{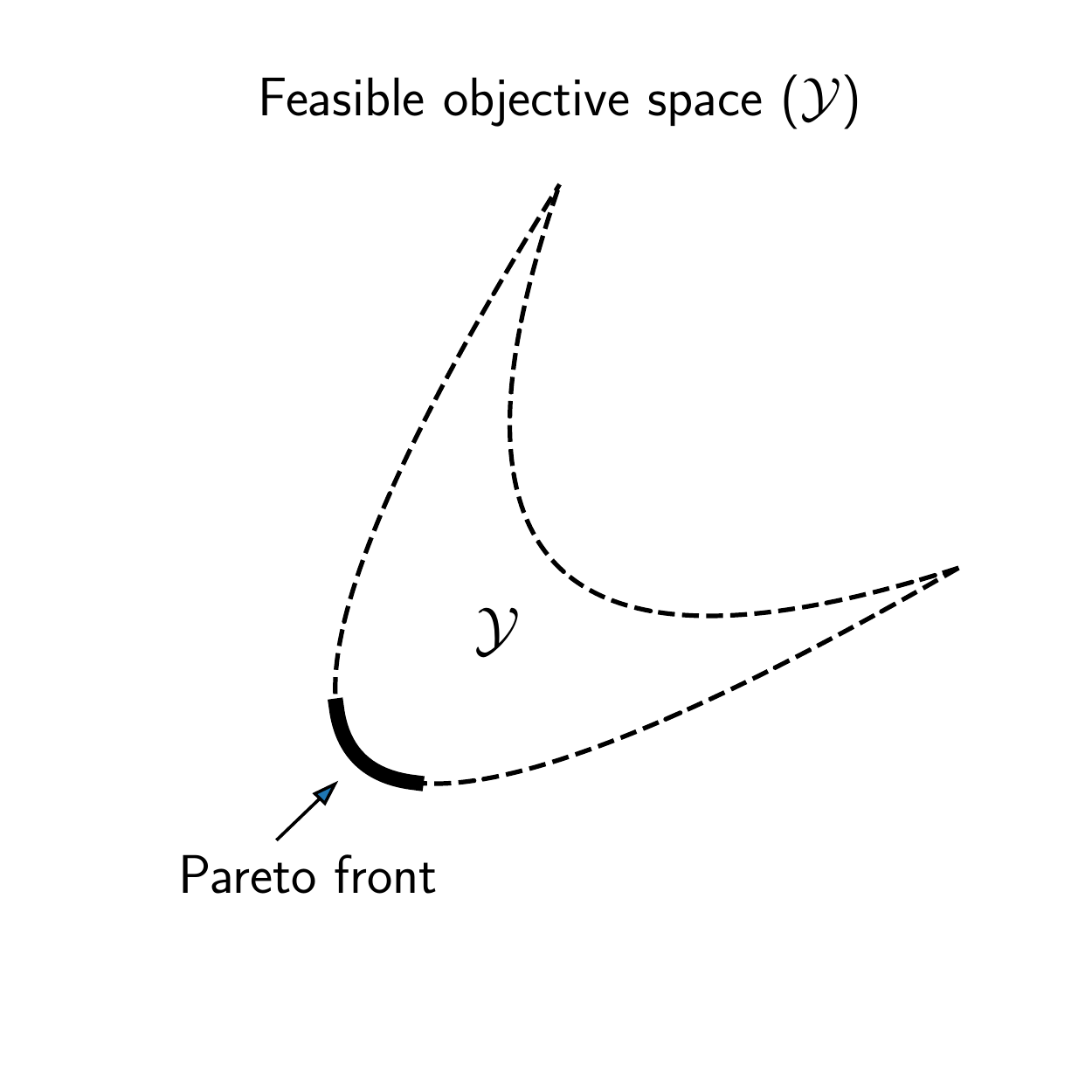}
\end{center}
\vskip -84pt 
\caption{The feasible design space ${\cX}$ (left) is mapped into an
intermediate simulation output space ${\cS}$ (center) via the black-box
simulation $\vS$, then to the feasible objective space ${\cY}$ (right) via
the algebraic objective function $\vF$.}
\label{fig:des-sim-obj}
\end{figure}

The next challenge that we need to address directly is Goal 2, which requires
making simulation outputs and function evaluations visible to the optimization
solver, for the problems defined in \Cref{eq:moso} and
\Cref{eq:parmoo}.
Our key insight here is to manage and model simulation outputs separately from
objective values and constraint penalties.  \Cref{fig:des-sim-obj} shows how
ParMOO treats MOOPs, with a black-box simulation function $\vS$ mapping into an
intermediate ``simulation output space'' ${\cS}$, before $\vF$ maps ${\cX}
\times {\cS}$ into the feasible objective space.

The abstraction depicted in \Cref{fig:des-sim-obj} allows us to handle the
problem formulations in \Cref{eq:moso} and \Cref{eq:parmoo}.  Instead of
applying surrogate modeling to the objectives using a dataset $\cD$ of design
point, objective value pairs, (e.g., fitting ${\hat \vF}_\cD(\vx_i) \approx
\vF(\vx_i)$) we model surrogates of the {\it intermediate simulation space}
$\vHatS(\vx_i) \approx \vS(\vx_i)$ using a dataset $\cD$ of design point, {\it
simulation output} pairs.  Then, similarly as in \Cref{eq:surrogate}, we use a solver to approximately minimize an acquisition function or
scalarization of the surrogate problem
\begin{equation}
\min_{\vx \in \cB} A_{\cD}(\vF(\vx, \vHatS(\vx)), \vHatSig(\vx))
\text{ subject to } \vG(\vx, \vHatS(\vx)) \leqq \vO
\label{eq:parmoo-constrained-iterate}
\end{equation}
where $\vHatSig$ is the uncertainty function for $\vHatS$ given $\cD$.

In each iteration, the {\tt SurrogateOptimizer} object has access to the
scalarized outputs $A_{\cD}(\vF(\vx, \vHatS(\vx)), \vHatSig(\vx))$, the
objective function $\vF$, and the simulation surrogate $\vHatS$, as well as the
constraint violation function $\vG$.  This information is sufficient to
implement a wide variety of composite-structure-exploiting optimization solvers.

\subsection{Managing Complexity through Problem Embeddings}

To address Goal 3, we need to handle problems involving continuous,
integer-valued, and categorical variables, as well as nonlinear constraints,
without requiring any problem-specific techniques or heavy modification to the
structure outlined in \Cref{sec:oop}.  To maintain usability and
maintainability (goal 4), we must be careful in the techniques that we use.

\begin{figure}[h!]
\begin{center}
\includegraphics[width=\textwidth]{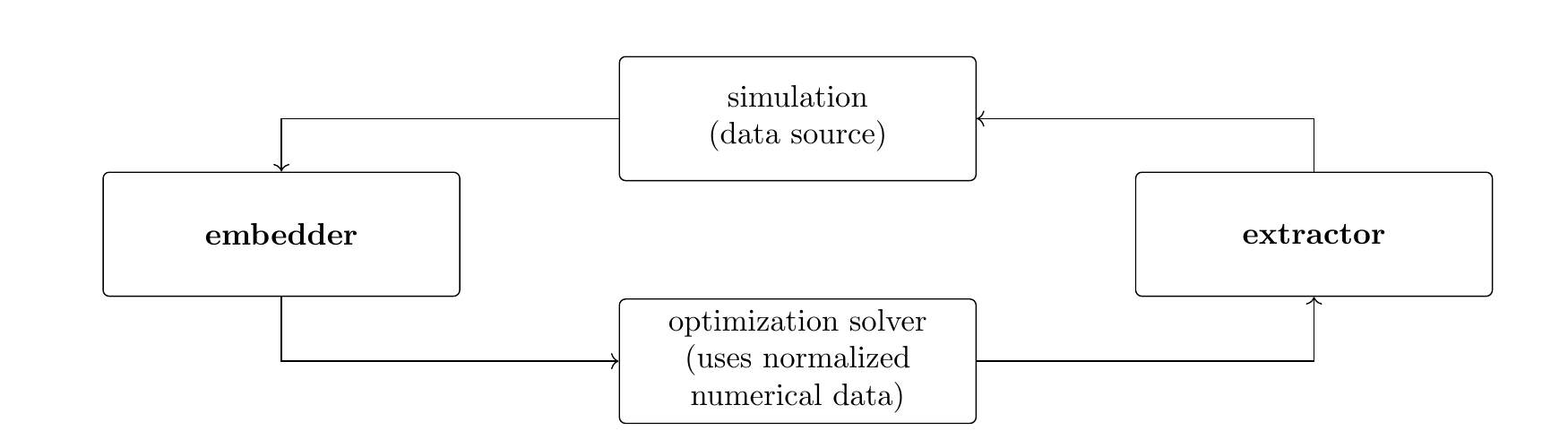}
\end{center}
\caption{
A depiction of the data flow in ParMOO.  The embedder module takes possibly
discrete values from the design space and maps them to a normalized latent
space $[0,1]^\ell$. The extractor module accepts candidate designs in the range
$[0,1]^\ell$ from the optimization module, and maps them back into $\cX$.
}
\label{fig:embedding}
\end{figure}

As discussed in \Cref{sec:embeddings-constraints}, for handling
mixed variables, our method of choice is to embed the problem into a continuous
latent space and then solve a relaxation of the problem.  ParMOO abstracts this
idea by creating a hidden embedding layer
based on built-in or user-provided embedder/extractor functions $\vEin$ and
$\vEout$ such that $\vEin(\vx) : \cX \rightarrow [0, 1]^\ell$ and
$\vEout(\vEin(\vx)) = \vx$;
see \Cref{fig:embedding}.  Upon input, all design variables are passed through
their respective embedders, $\vEin$, and any variables taking on discrete values
are relaxed via the surrogate model.
After ParMOO solves the surrogate optimization problem, all candidate design
points are extracted from the latent space back into the
original design space via $\vEout$
and binned to their nearest legal values if necessary.  ParMOO provides a
default embedding for categorical variables, which utilizes a one-hot encoding,
followed by a dimension reduction procedure to eliminate unused dimensions of
the latent space (i.e., two different latent variables representing the
same categorical variable cannot both have nonzero values at the same time).
However, this technique is not scalable to a large number of categorical
variables.  Therefore, for problems with large numbers of categorical
variables, ParMOO also allows users to provide a custom embedding procedure.
For continuous and integer design variables, the default embedding is based
upon a simple rescaling to the range $[0, 1]$.

The primary advantage of this approach is that it maintains a level of
simplicity in the API and allows us to focus on continuous optimization
techniques
regardless of the problem type.  In particular, there is no need to use
customized search, surrogate modeling, or optimization techniques for solving
mixed-variable problems (although doing so may improve performance) since the
default techniques are designed to operate on the rescaled latent space.  The
flexibility to accept a custom embedding tool is also useful for domains such
as molecular engineering, where these tools are available \cite{moriwaki2018}.
As previously discussed, for problems involving large numbers of
categorical variables without any domain-specific embedding tools, the default
method may become inefficient and it may be advantageous to take a completely
different approach that is not based on continuous MOSO algorithms.

In order to incorporate a $p$-dimensional nonlinear constraint function $\vG$,
we apply a cumulative penalty function,
$P(\vx) = \lambda\sum_{i=1}^p \max(G_i(\vx, \vS(\vx)), 0)$ 
or $P(\vx) = \lambda\sum_{i=1}^p \max(G_i(\vx, \vHatS(\vx)), 0)$, to the
surrogate problem.  This penalty is
added to {\it all} of the objective functions when solving
\Cref{eq:parmoo-constrained-iterate}.
Combining with the embedder/extractor functions described above, we can
reduce a nonlinearly constrained, partially discrete feasible design space $\cX$
to a bound-constrained, well-conditioned, continuous MOSO problem
\begin{equation}
\min_{\vx \in [0,1]^\ell} A_{\cD}(\vF(\vEout(\vx), \vHatS(\vx)) +
\lambda\sum\max(\vG(\vEout(\vx), \vHatS(\vx)), \vO), \vHatSig(\vx)),
\label{eq:parmoo-iterate}
\end{equation}
where the $\max$ is taken elementwise.  The penalty parameter $\lambda>0$ is
progressively increased, on an exponential schedule, whenever the solution to
\Cref{eq:parmoo-iterate} is not expected to be feasible.  This is similar to
the progressive barrier approach described in \cite{audet2009}, and has been
successfully applied with simulation-based constraints in \cite{neveu2023}.
This approach has been shown to be effective when the constraints are
hard to satisfy \cite{chen2023}.

One of the downsides of this approach, is that the penalty function uses
constraint violations, which implicitly requires the constraint functions in
$\vG$ to be both relaxable and quantifiable.  Therefore, in ParMOO, only the
upper/lower bound constraints on each design variables can be unrelaxable,
which is a necessary compromise to achieve reasonable performance on heavily
constrained problems.

\subsection{Flexibility and Extensibility to Novel Workflows}

Goal 4 was to make our framework flexible enough that it is easy to deploy in a
variety of scientific workflows.  Key examples include HPC and wet-lab
environments.

To achieve this goal, we have designed ParMOO  to issue simulation evaluations
exclusively via the wrapper function {\tt MOOP.evaluateSimulation()}, which is
exclusively called from within the solver routine {\tt MOOP.solve()}.  In
situations where only the simulation command must be changed for integration
with existing workflow technology, this can be achieved by extending the {\tt
MOOP} class and overwriting the {\tt MOOP.evaluateSimulation()} method.  In
other situations, where control over the frequency with which simulation
evaluations are distributed is required (for example, when batching simulation
evaluations), the entire {\tt MOOP.solve()} method can be overwritten.
The latter approach is taken in ParMOO's {\tt libE\_MOOP} class, which extends
the {\tt MOOP} class and overwrites the {\tt solve} method to dynamically
distribute simulation evaluations on HPC systems using the libEnsemble library
\cite{hudson2022a,hudson2022b}.  The {\tt libE\_MOOP} class is currently
the recommended method for achieving scalable parallelism with ParMOO.

\subsection{The Resulting Framework}

Putting everything together, we present our framework for solving MOSO problems
as implemented in ParMOO.
In each iteration, ParMOO uses \Cref{alg:parmoo-iterate} to generate the next
batch of candidates (implemented in the {\tt MOOP.iterate()} method).

\begin{algorithm}[h!]
    \SetKwInOut{Input}{input\hskip 0.7em}
    \SetKwInOut{Output}{output}

\IncMargin 1em

    \Input{$k \geq 0$ is the current iterate.}
    \Input{$q_0 \geq 0$ is an initial search budget.}
    \Input{$q$ is the batch size, i.e., the number of
           {\tt AcquisitionFunction}s.}
    \Input{$n$ is the dimension of the design space and $\ell$ is the dimension
           of the latent space.}
    \Input{$\vEin : \cX \rightarrow [0, 1]^\ell$ is the cumulative embedder
           function for all design variables.}
    \Input{$\vEout : [0, 1]^\ell \rightarrow \cX$ is the extractor function
           (inverse of $\vEin$).}
    \Input{$\vS : \cX \rightarrow \R^m$ is the simulation function.}
    \Input{$\vF : \cX \times \R^m \rightarrow \R^o$ is the objective function.}
    \Input{$\vG : \cX \times \R^m \rightarrow \R^p$ is the constraint function.}
    \Input{$\cD^{(k)}$ is the set of all $(\vx, \vS(\vx))$ pairs
           evaluated prior to iteration $k$.}
    \Input{{\tt GlobalSearch.search}($q_0$, $\ell$) is a procedure that
           generates a design-of-experiments or sample of size $q_0$ in $[0,
           1]^\ell$.}
    \Input{{\tt SurrogateFunction.fit}($\cD$) fits the surrogates
           $\vHatS : [0, 1]^\ell \rightarrow \R^m$, and (if applicable) their
           uncertainty function
           $\vHatSig : [0, 1]^\ell \rightarrow \R^m$.}
    \Input{{\tt SurrogateFunction.improve}($\cT$) produces a model
           improving step for $\vHatS$ in $\cT \subset [0, 1]^\ell$ (if
           applicable).}
    \Input{\{$A_{\cD}^{(i)}\}_{i=1}^q$ is a set of $q$
           {\tt AcquisitionFunction}s.}
    \Input{{\tt SurrogateOptimizer.solve}($A_\cD$, $\vF$, $\vG$, $\vHatS$,
           $\vHatSig$, $\vx^{(0)})$) solves \Cref{eq:parmoo-iterate},
           iterates toward the solution to \Cref{eq:parmoo-iterate} starting
           from $\vx^{(0)}$, or calls the
           {\tt SurrogateFunction.improve($\cT$)} method and provides a LTR
           $\cT$.}
    \Output{The batch $\cC^{(k)}$ of candidate design points after iteration
            $k$.}

    $\cC^{(k)} \leftarrow \emptyset$\;
    \eIf{$k = 0$}{
        $\cC' \leftarrow ${\tt GlobalSearch.search}$(q_0, \ell)$\;
        \ForEach{$\vy \in \cC'$}{
            $\cC^{(0)} \leftarrow \cC^{(0)} \cup \{\vEout(\vy)\}$\;
        }
    }{
        $\cD' \leftarrow \{(\vEin(\vx), \vS(\vx)) : \text{ for all }
            (\vx, \vS(\vx)) \in \cD^{(k)}\}$\;
        $\vHatS, \vHatSig \leftarrow$ {\tt SurrogateFunction.fit}($\cD'$)\;
        \For{$i \leftarrow 1$ \KwTo $q$}{
            $\vx^{(0,i)} \leftarrow \arg\min_{(\vx,\vS(\vx))\in\cD'}$
                $A_{\cD'}(\vF(\vx, \vS(\vx)), \vO))$\;
            $\vy^{(k,i)} \leftarrow$ {\tt SurrogateOptimizer.solve}(
                $A_{\cD'}$, $\vF \circ \vEin$, $\vG \circ \vEin$, $\vHatS$,
                $\vHatSig$, $\vx^{(0,i)}$)\;
            \eIf{$\vy^{(k,i)} \not\in \cC^{(k)}$}{
                $\cC^{(k)} \leftarrow \cC^{(k)} \cup \{\vEout(\vy^{(k,i)})\}$\;
            }{
                $\vy' \leftarrow $
                {\tt SurrogateFunction.improve}($[0,1]^\ell$)\;
                $\cC^{(k)} \leftarrow \cC^{(k)} \cup \{\vEout(\vy')\}$\;
            }
        }
    }
    \Return $\cC^{(k)}$\;

\DecMargin 1em
    \caption{A single iteration of ParMOO's {\tt MOOP.iterate()} method
    \label{alg:parmoo-iterate}
    }
\end{algorithm}

To fill in the various components of \Cref{alg:parmoo-iterate}, the user must
create and fully populate a {\tt MOOP} object.  When initializing the {\tt
MOOP}, the user must specify a {\tt SurrogateOptimizer} for suggesting
candidates for solving or iterating upon \Cref{eq:parmoo-iterate} and a
dictionary of hyperparameters.
Next, the user must add problem details, including 
\begin{itemize}
\item $n$ design variable dictionaries, each specifying a design
variable for the problem and
(if applicable) a custom {\tt Embedder} and {\tt Extractor}
routine;
\item $s$ simulation dictionaries, each defining a simulation/data source
for the problem, plus the number of intermediate outputs ($m_i$), the
{\tt GlobalSearch}
technique used to sample that simulation's output space, and the
{\tt SurrogateFunction}
used to model that simulation's output
and either evaluate uncertainties or perform model improvement iterates;
\item $o$ objective dictionaries, each specifying an algebraic function that
can be used to calculate the objective values from the design variables and
simulation outputs and (optionally) the gradient of that objective;
\item $p$ constraint dictionaries, identical to the objective dictionaries
but for the purpose of evaluating relaxable constraint violations
(i.e., the $\vG$ functions); and
\item $q$ acquisition function dictionaries, each specifying an
{\tt AcquisitionFunction}
that can be used to scalarize the problem, where $q$ also determines the batch
size for parallel evaluations.
\end{itemize}

\begin{figure}[ht]
\begin{center}
\includegraphics[width=0.99\textwidth]{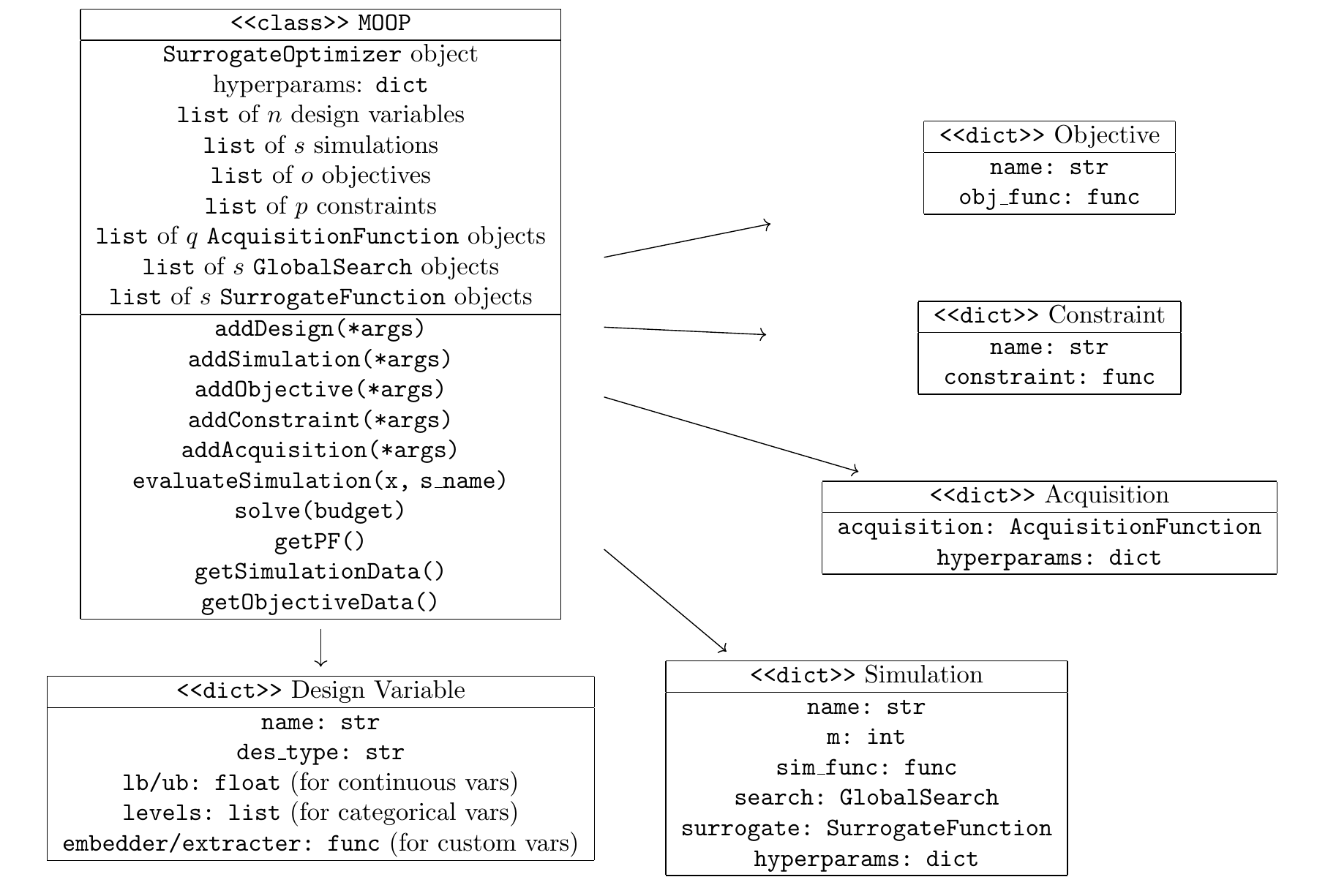}
\end{center}
\caption{UML diagram outlining the key dictionaries, components,
and methods that make up a {\tt MOOP} object and its contents.}
\label{fig:uml2}
\end{figure}

A UML diagram outlining this framework is shown in \Cref{fig:uml2}.  Note that
only the relevant public methods discussed in this section are shown.  In the
{\tt MOOP} class's implementation, several additional public and private helper
methods are included, as well as additional ``setter'' methods to facilitate
saving, loading, logging, and checkpointing \cite{chang2023a}.

\subsection{Limitations of the Framework}
\label{sec:limitations}

We have endeavored to cover a great number of MOSO algorithms and techniques
and make them easily accessible for a wide variety of applications.  While we
have succeeded for a wide breadth of algorithms, our framework is not without
limitations.  For MOSO algorithms that do not use acquisition functions or
surrogates, as discussed in \Cref{sec:other-techniques}, it would likely be
inefficient to implement these in ParMOO's framework and thus shoehorn their
iterations into \Cref{alg:parmoo-iterate}.  Additionally, as acknowledged in
\Cref{sec:embeddings-constraints}, ParMOO's framework is fundamentally based
upon techniques for continuous, bound-constrained MOSO and handles mixed
variables and nonlinear constraints by reducing these problems to
bound-constrained, continuous MOSO problems.  For problems with no continuous
variables or unrelaxable constraints, this is likely not the most effective
approach and other frameworks and classes of algorithms may be better suited.

Although MOSO is typically a derivative-free endeavor, we have
made efforts to offer limited support for partial derivative information, since
the algebraic functions $\vF$ and $\vG$ could be differentiable and ParMOO's
optimization solvers can use their gradient information when available.
However, this can become misleading in the sense that ParMOO is designed to
solve problems involving at least {\it some} computationally expensive
black-box simulations.  If there is no black-box simulation function involved
in the problem definition and if gradients are available for all components of the
objective and constraint definitions, then other approaches not covered here
would likely be more effective.

Finally, it is worth noting that ParMOO is fundamentally a framework, not
an algorithm.  Therefore, much of the work in exploiting problem structures and
implementing quality algorithms is left to the user.  While ParMOO's framework
and development is mature, the variety of techniques available is still under
active development and not all of the techniques listed in
\Cref{sec:techniques} are available at this time.
Careful analysis and clever algorithms are needed to gain a theoretical
advantage when exploiting most problem structures, although we will show in
\Cref{sec:fayans} and \Cref{sec:chem} that ParMOO can still gain some clear
advantage in most applications with relatively little work.

\section{Performance Summary}
\label{sec:dtlz2}
In this section we apply a trust-region solver built in ParMOO to a common
academic test problem from the literature in order to assess its parallel
scaling.  Before doing so, we would like to caution about the limitations of
this study.  First, as described in \Cref{sec:framework}, ParMOO is a framework
and software library for implementing solvers, but not an algorithm in and of
itself; there can be significant variations in ParMOO's performance
depending on the methods used.  Second, the focus of this paper and ParMOO's
design is on achieving performance on real MOSO applications.  Since the
academic problem used here does not necessarily reflect a real-world MOSO
application, performance on this problem should be understood as validating
that ParMOO achieves reasonable performance on solvable problems, and not that the
particular ParMOO solver described is ``better'' or ``worse'' than existing
methods.  Thus, we have intentionally chosen a single relatively easy problem
from the multiobjective optimization literature, the DTLZ2 problem from the widely-used DTLZ
test suite \cite{deb2002b}.  We use the variation of this problem with 3
objectives and 10 design variables.  The problem has no constraints, no mixed
variables, and no composite structures of the form shown in
\Cref{fig:des-sim-obj} to exploit.  The simplicity of the problem
allows us to focus on the correctness of our implementation and parallel
scaling under generic conditions.

Since the entire DTLZ2 problem is algebraic, we must create an artificial
simulation.  Therefore, we have created a wrapper for the DTLZ2 function that
turns it into a black-box simulation by creating an artificially long runtime.
This is done by uniformly sampling a runtime between one and three seconds,
then waiting for that amount of time using Python's built-in {\tt time.sleep()}
function.  The objectives ($\vF$-functions) are then identity mappings from
this artificial simulation's three output fields.

As our sample solver, we have created a parallel ParMOO solver for this problem
by instantiating the {\tt libE\_MOOP} builder class.  The solver has been
populated using the {\tt TR\_LBFGSB} implementation of the {\tt
SurrogateOptimizer} to solve a trust-region-constrained surrogate problem with
L-BFGS-B; a combination of the randomized epsilon-constraint scalarization {\tt
RandomConstraint} and fixed-weight weighted sum scalarization {\tt FixedWeight}
for our {\tt AcquisitionFunction}s; a {\tt GaussRBF} (Gaussian RBF) {\tt
SurrogateFunction}; and a Latin hypercube {\tt GlobalSearch}.  By varying the
number of {\tt AcquisitionFunction} objects, we are able to create variations
of this ParMOO solver that will generate batch sizes of 8, 16, and 32 candidate
simulations per iteration.  Since ParMOO parallelizes simulation evaluations,
this batch size determines the largest amount of parallelism available in a
single iteration.

Since the Latin hypercube sampling and default model-improvement procedures are
randomized, we have performed 5 runs of the above-defined ParMOO solver with
different random seeds for each run.  We have solved with each random seed and
each batch size serially and with 2-, 4-, and 8-way parallelism during simulation
evaluation with a budget of 1,000 total simulation evaluations.

To validate our results, we have calculated the widely-used hypervolume
performance indicator at the end of each run \cite{audet2020}, averaged these
values, and compared across runs.  For the hypervolume indicator, larger values
are better, with the maximal value being the total volume between the Pareto
and a reference point (in this case, the point $(1, 1, 1)$).  Hypervolume
values are notoriously difficult to interpret, so for comparison, we have also
solved DTLZ2 on the same budget using the widely-used NSGA-II implementation in
pymoo \cite{deb2002a,blank2020}.

The results of the comparison are shown in \Cref{tab:dtlz2-hv}.  Even without
exploiting any structures, this shows that ParMOO solver's performance on DTLZ2 (according
to the hypervolume indicator) is consistently higher than pymoo after 1,000 simulation evaluations, indicating that ParMOO achieves at least
reasonable performance on academic benchmark problems.  There is slight
variation between the performance of the ParMOO solvers of different batch
sizes, however, we would recommend against drawing conclusions on such limited
data, especially since DTLZ2 is not a particularly representative application
for ParMOO.

\begin{table}[h]
\begin{center}
\begin{tabular}{c|cccc}
Method & pymoo & ParMOO-8 & ParMOO-16 & ParMOO-32\\
\hline
Hypervolume & 0.28	&	0.33	&	0.33	&	0.37\\
\end{tabular}
\vskip 12pt
\end{center}
\caption{Hypervolume indicator for pymoo and ParMOO with batch sizes 8, 16, 32.
Larger values are better.}
\label{tab:dtlz2-hv}
\end{table}

\begin{figure}[h!]
\begin{center}
\includegraphics[width=0.75\textwidth]{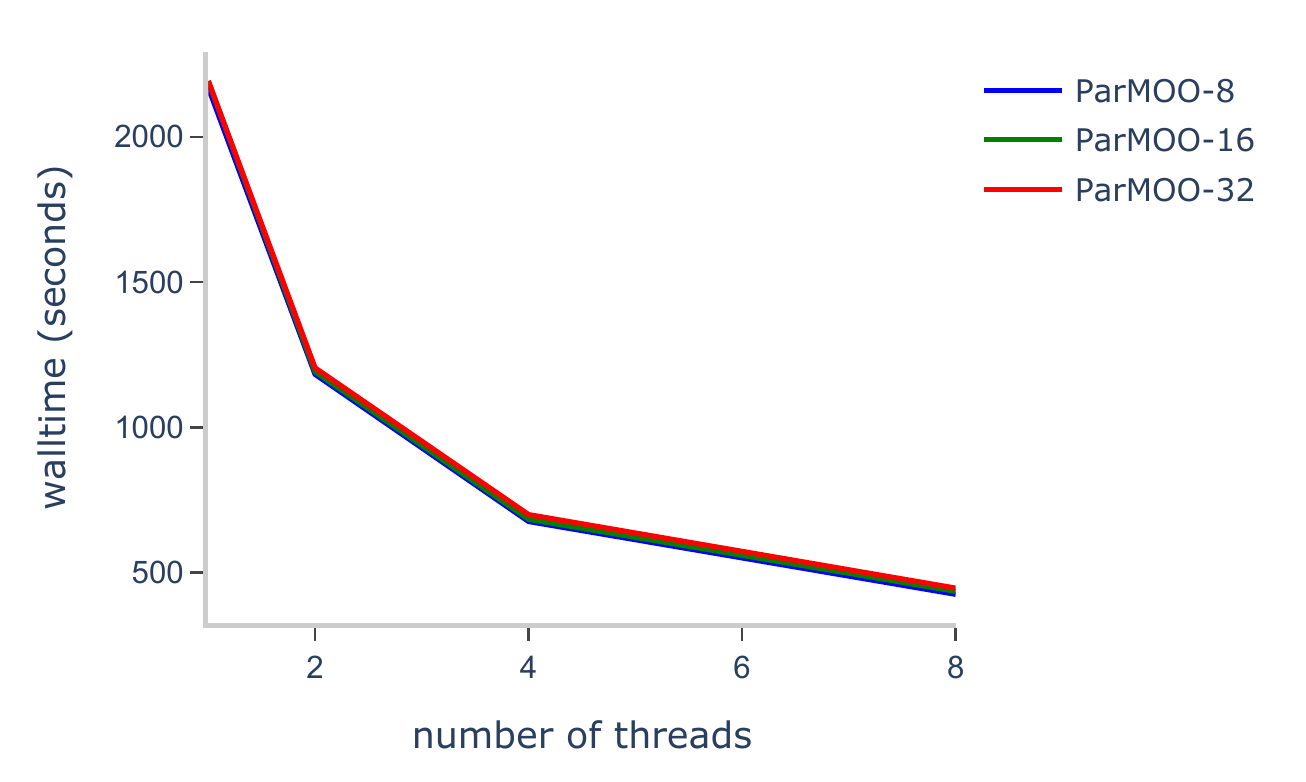}
\end{center}
\caption{Walltimes when performing 1,000 1-3 second simulation evaluations in
ParMOO (with batch sizes 8, 16, and 32) with increasing number of threads.}
\label{fig:dtlz2-walltime}
\end{figure}

ParMOO does not change \Cref{alg:parmoo-iterate} when running in parallel, it
only parallelizes the simulation evaluations in a batch.  Therefore, we have
not observed any changes in its convergence (as measured by the hypervolume
indicator) when increasing the number of threads other than minor and
uncorrelated fluctuations due to random number generation.  Therefore, if we
are able to reduce the walltime to perform a 1,000 simulation run with any of
the above solvers, this marks a distinct parallel advantage.
In \Cref{fig:dtlz2-walltime}, we consider the total walltimes when running each
ParMOO batch size with increasing number of threads.  Note that since the
average simulation time is 2 seconds, the expected total walltime when
performing 1,000 total simulation evaluations serially would be approximately
2,000 seconds plus any iteration costs incurred by ParMOO.

As seen in \Cref{fig:dtlz2-walltime}, ParMOO iteration computations incur a bit
of overhead for this solver configuration and problem size (note that the
single threaded walltime is a bit over 2,000 seconds) but the walltime to
solution decreases proportionally with the number of threads.  This parallelism
is slightly less than perfect strong scaling and is limited by the
candidate batch size.  Still, for most reasonable solver configurations, when
performing 8 simulations in parallel requiring an average 2 seconds per
simulation ($\pm$1 second variation), one could expect to perform 1,000
simulations in ParMOO in just under 500 seconds.

\section{Case Study: Solving a Multiobjective Inverse Problem}
\label{sec:fayans}
In this section we apply ParMOO to a multiobjective inverse problem, where the
goal is to tune the parameters of the Fayans energy density functional (EDF)
model based on experimental data.  This problem is fully described in
\cite{bollapragada2020}.  The forward model is expensive to evaluate and not
publicly available.  Therefore,  to maintain reproducibility of results, we
optimize a synthetic problem based on a neural network model that was trained
on the dataset above.

\subsection{Background on Fayans EDF Calibration}

We now review important aspects of the Fayans EDF calibration.  Let $T_\vx : s
\rightarrow t$ denote the Fayans EDF model.  A single forward evaluation of
$T_\vx$ accepts a normalized parameter vector $\vx\in\R^{13}$ and an
input $s$ and produces a physical observation $t$.  Given $198$ observational
data pairs $(s_1, t_1),\ldots,(s_{198}, t_{198})$, the goal of the Fayans EDF
model calibration is to find ``good'' parameter vectors $\vx^\star \in
\R^{13}$ such that $T_{\vx^\star}(s_i) \approx t_i$ relative to a given
standard error $\sigma_i>0$ for $i=1, \ldots, 198$.

In traditional Fayans EDF calibration, it is assumed that all observations have
been normalized by the standard errors $\bm{\sigma}$ such that they have equal
and independent normalized errors.  It is then reasonable to minimize the
$\chi^2$ loss across all 198 observations via the single-objective formulation 
\begin{equation}
\min_{\vx\in\R^{13}} \sum_{i=1}^{198} \left(\frac{T_\vx(s_i) - t_i}{\sigma_i}\right)^2.
\label{eq:fayans-single}
\end{equation}

In \cite{bollapragada2020}, however, it is noted that the contributions to the
$\chi^2$ loss may vary across 9 different observational types.  Since the
9-objective formulation is prohibitively expensive to solve, we consider a
simplified three-objective formulation, where related observational types are
combined, leaving just three observational classes.  This allows us to define
the following multiobjective formulation of the Fayans EDF calibration:
\begin{equation}
\min_{\vx\in\R^{13}} F_{\vx,j}, \qquad j=1,2,3,
\label{eq:fayans-multi}
\end{equation}
where $F_{\vx,j} = \sum_{i\in\Phi_j}\left(\frac{T_\vx(s_i) -
t_i}{\sigma_i}\right)^2$ and $\Phi_1, \Phi_2, \Phi_3$ is a partitioning of
$\{1, \ldots, 198\}$ based on the three observational classes described above.

In this section we solve the multiobjective formulation of the Fayans EDF
calibration problem \cref{eq:fayans-multi}.  In this problem, evaluation of the
forward model for all $198$ observations $[T_\vx(s_1), \ldots, T_\vx(s_{198})]$
is viewed as a single computationally expensive black box.  The three loss
functions $F_{\vx,j}$ ($j=1,2,3$) all have
a composite structure, specifically, a sum-of-squares structure.

\subsection{The Neural Network Residual Model}

For our experiments, since $T_\vx$ is expensive and we strive for
reproducibility, we fit a multilayer perceptron (MLP) to approximate the
standardized residual functions $R_i(\vx) = \frac{T_\vx(s_i) - t_i}{\sigma_i}$,
$i=1,\ldots,198$.  We do this using a dataset of labeled observations for the
residual functions $R_1(\vx), \ldots, R_{198}(\vx)$ for 52,079 distinct values
of $\vx$.  This dataset was gathered by \cite{bollapragada2020} and, as
described in that paper, comes from running multiple starting design points
with several different single-objective solvers on the single-objective
formulation of the problem \cref{eq:fayans-single}.  Therefore, we note that
the dataset is not uniformly distributed and has increased density in the
neighborhood of several local solutions to the single-objective problem
\cref{eq:fayans-single}.

The primary challenge when training the MLP used in this section is ensuring
that the prediction accuracy is high for optimal values of the multiobjective
problem \cref{eq:fayans-multi}.  In particular, large regions of the parameter
space  were determined to be unstable for the single-objective formulation, and
therefore no observational data is available.  Additionally, within the stable
region, for several values of $i$, the observational data for $R_i(\vx)$ could
range in magnitude from less than $1$ to over $10^{20}$.  While prediction
errors that are greater than $1$ in magnitude would be acceptable (and
expected) for these large residuals, it is essential for our fidelity to the
original problem that we maintain prediction errors on the order of $10^{-2}$
in the neighborhood of the true solutions.

To handle infeasible regions of the parameter space where no observational data
is available, we impose bound constraints for ParMOO, following the constraints
given by \cite{bollapragada2020}, shown in \Cref{tab:fayans-bounds}.  To handle
the wide range on observational values, we apply a double-logarithmic
transformation to large values of $R_i(\vx)$, followed by a $\tanh$
transformation.  Because this collapses huge error values, this transformation
ensures that, during training, the MLP places much higher importance on
matching observations where the residual function $R_i(\vx)$ has a low score.

\begin{table}[h]
\begin{center}
\begin{tabular}{c|cc}
Variable Name & Lower Bound & Upper Bound\\
\hline
$\rho_{eq}$     & 0.146     & 0.167     \\
$E/A$           & -16.21    & -15.50    \\
$K$             & 137.2     & 234.4     \\
$J$             & 19.5      & 37.0      \\
$L$             & 2.20      & 69.6      \\
$h^v_{2-}$      & 0.0       & 100.0     \\
$a^s_+$         & 0.418     & 0.706     \\
$h^s_{\nabla}$  & 0.0       & 0.516     \\
$\kappa$        & 0.076     & 0.216     \\
$\kappa'$       & -0.892    & 0.982     \\
$f^{\xi}_{ex}$  & -4.62     & -4.38     \\
$h^\xi_{+}$     & 3.94      & 4.27      \\
$h^\xi_{\nabla}$ & -0.96    & 3.66      \\
\end{tabular}
\vskip 12pt
\end{center}
\caption{Upper and lower bounds on the stable region for the parameter space
for the Fayans EDF model $T_\vx(s)$, as reported by \cite{bollapragada2020}.
We  have access to observational data only in these ranges, and therefore the
trained MLP will  be accurate only within these bounds.}
\label{tab:fayans-bounds}
\end{table}

For the network architecture, we define an MLP with $13$ inputs, $198$ outputs,
and 2 hidden layers, with $256$ nodes per layer.  We then apply $\tanh$
activation functions for every layer of the network (including the output
layer) and train the network using the transformed data defined above.  For
validation, we first stratify the complete observational database by low/high
residuals across all three observational classes defined in
\cref{eq:fayans-multi} and then withhold 5\% of the data across all
stratifications for our validation set.  The network was trained in {\tt keras}
\cite{chollet2015} using 5,000 epochs of RMS-prop.  We verify that the trained
model obtains low relative error across all stratified residual ranges in the
validation set.  Most importantly, after descaling the outputs to their
original ranges, on the bottom stratification where the residuals are lowest in
magnitude across all three observational classes, we obtain a mean absolute
error (MAE) of just 0.036 on the validation set, which is acceptable accuracy
for this problem.

In the remainder of this section we use the trained {\tt keras} model described
above as a synthetic representation of the original Fayans EDF calibration
problem.

\subsection{Exploiting Structure in ParMOO}

The sum-of-squares structure is well studied in the single-objective black-box
optimization literature.
In typical software implementations of this strategy \cite{wild2017}, special
care is taken to ensure that all geometric conditions on the interpolation
nodes hold, thus requiring the algorithm to occasionally take
geometry-improving evaluations and using a sequence of local (trust-region)
approximations.
However, it has also been shown that simply modeling this composite structure
typically leads to fewer simulation evaluations in practice
\cite{astudillo2021,khan2018,larson2023}.

In ParMOO's built-in {\tt LocalGaussRBF} surrogate model implementation, the
trust-region radius is adaptively chosen based on the distance to the $(n+1)$th
nearest neighbor of the current iterate.  Then, whenever the solution to the
trust-region-constrained surrogate problem fails to produce a sufficient
decrease in the scalarized objective value,
the model-improving step is triggered, which draws a random sample from a
distribution whose variance is highest in directions of low variance in the
current simulation dataset.  This approach does not guarantee any asymptotic
convergence rate advantage (indeed, such a convergence rate would be difficult
to even define in the multiobjective setting).
We will see empirically, however, that simply modeling $\vR(\vx)$ separately
from $\vF(\vx)$ greatly accelerates the practical convergence of our
structure-exploiting multiobjective solver.

\subsection{The ParMOO Fayans EDF Solvers}

To demonstrate its ability to utilize parallel function evaluations (as we
would need to if we were evaluating the true Fayans EDF model instead of our
synthetic MLP model), we use ParMOO's {\tt libEnsemble} interface to distribute
simulation evaluations with a batch size of 10.
This is achieved by using $9$ epsilon-constraint acquisition functions that
target the solutions to \Cref{eq:fayans-multi} and one fixed weighting that
targets the solution to the single-objective formulation in
\Cref{eq:fayans-single}, since it is equivalent to the $\chi^2$ loss.
To define the problem, we create a {\tt libE\_MOOP} object in ParMOO and add
$13$ continuous design variables, with names and bounds as given in
\Cref{tab:fayans-bounds}.

To define a structure-exploiting variant as described in the preceding section,
we then add a single simulation function that models all $198$ components of
our synthetic MLP residual function using LTR constrained
Gaussian radial basis function (RBF) models and initializes its observational
database using a $2,000$-point Latin hypercube sample (LHS).  We next add three
differentiable algebraic objective functions, which calculate the sum of
squares from the simulation outputs across each of the three different
observational classes.  Then, to allow ParMOO to focus on interesting regions
of the design space where none of the observational classes result in terrible
performance, we add three nonlinear constraints, each enforcing that the
squared residuals in a particular observable class should not exceed 10$\times$
its optimal values reported in \cite{bollapragada2020}.  ParMOO then solves the
surrogate optimization problems within the LTR using the L-BFGS-B
\cite{zhu1997} implementation in SciPy \cite{virtanen2020}.

To compare with the performance when the structure is not exploited, as would
be the case with other existing multiobjective software, we define a
``black-box'' variation of this problem.  For the black-box variant, the
$198$-output simulation described above was replaced by a $3$-output
``black-box'' simulation, where the sum of squares has already been computed
across all three observational classes.  Then, for each objective, we provide
an identity map from each simulation output to each objective.  Otherwise, we
define the problem equivalently as with the structured variant.

In both cases, ParMOO is then run for 800 iterations.  Since we begin with a
$2,000$-point LHS search and provide a batch size of 10 points per iteration,
this results in a total budget of 10,000 simulation evaluations for each
solver.  We note that for a typical simulation optimization problem, this
budget of 10,000 simulation evaluations is unrealistically large.  In fact,
because of the complexity of our Gaussian RBF surrogate models and cost of
solving the surrogate optimization problem, for this large a budget our
iteration costs are extremely high.  However, it is worth running to such a
large simulation budget in order to understand the solver's performance in the
limit.

Note that because of the complexity of solving the surrogate optimization
problem and the nonnegligible cost of evaluating the {\tt keras} model, such a
solve requires a substantial amount of compute time.  To account for randomness
in the LHS searches and weight initializations, we perform 5 runs of each
solver and average the performance.

For reproducibility, this entire experimental setup, including the trained MLP
that was used as a synthetic problem representation, is available at
\url{https://github.com/parmoo/parmoo-solver-farm/tree/main/fayans-model-calibration-2022}.
Note that for ease of use and compatibility reasons, we have transferred our
{\tt keras} model's weights into an equivalent {\tt torch} model
\cite{paszke2019}, and this is reflected in the GitHub repository.

\subsection{Results}
\label{sec:fayans-results}

\begin{figure}[h!]
\begin{center}
\includegraphics[width=0.75\textwidth]{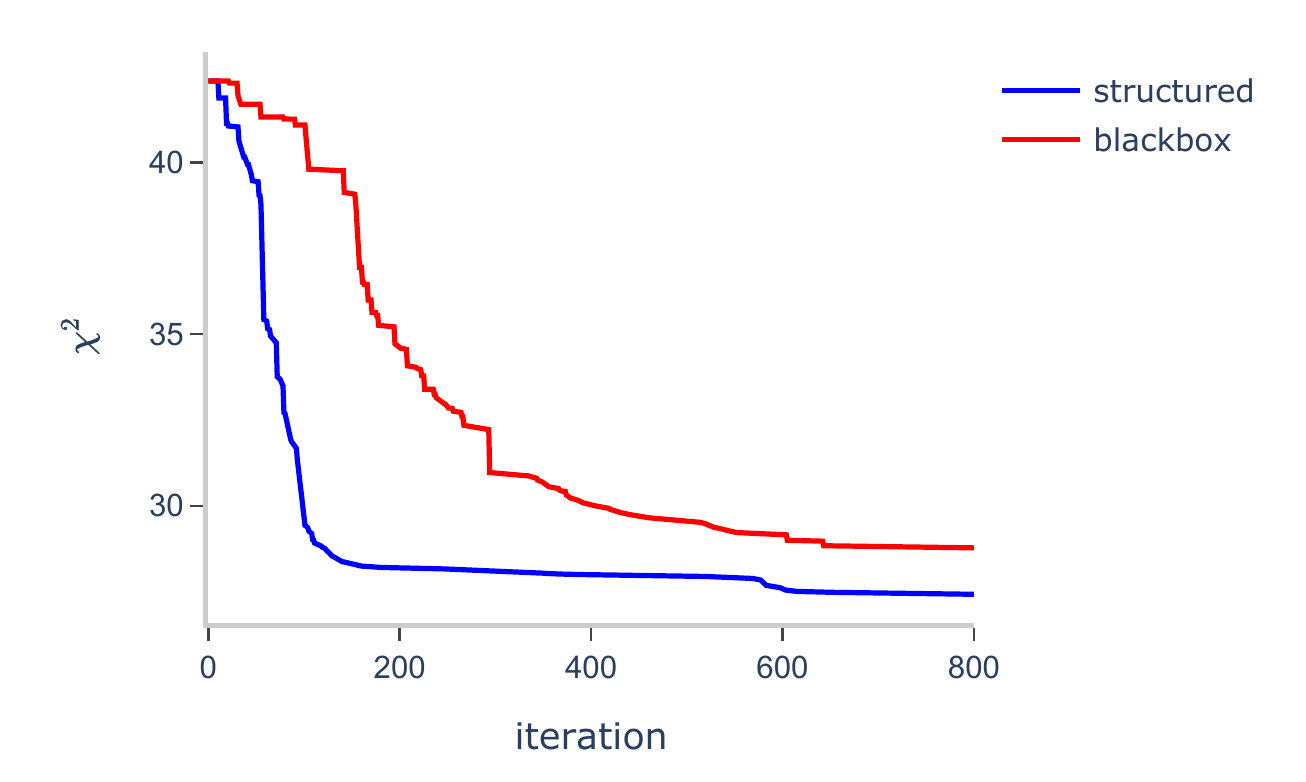}
\end{center}
\caption{Iteration vs.\ $\chi^2$ loss when solving the Fayans EDF calibration
with ParMOO, exploiting the sum-of-squares structure (structured) and with
a standard (black-box) approach.  The total simulations used by the end of
iteration $k$ are calculated as $2000 + 10k$.
}
\label{fig:fayans-results-chi2}
\end{figure}

\begin{figure}[h!]
\begin{center}
\includegraphics[width=0.75\textwidth]{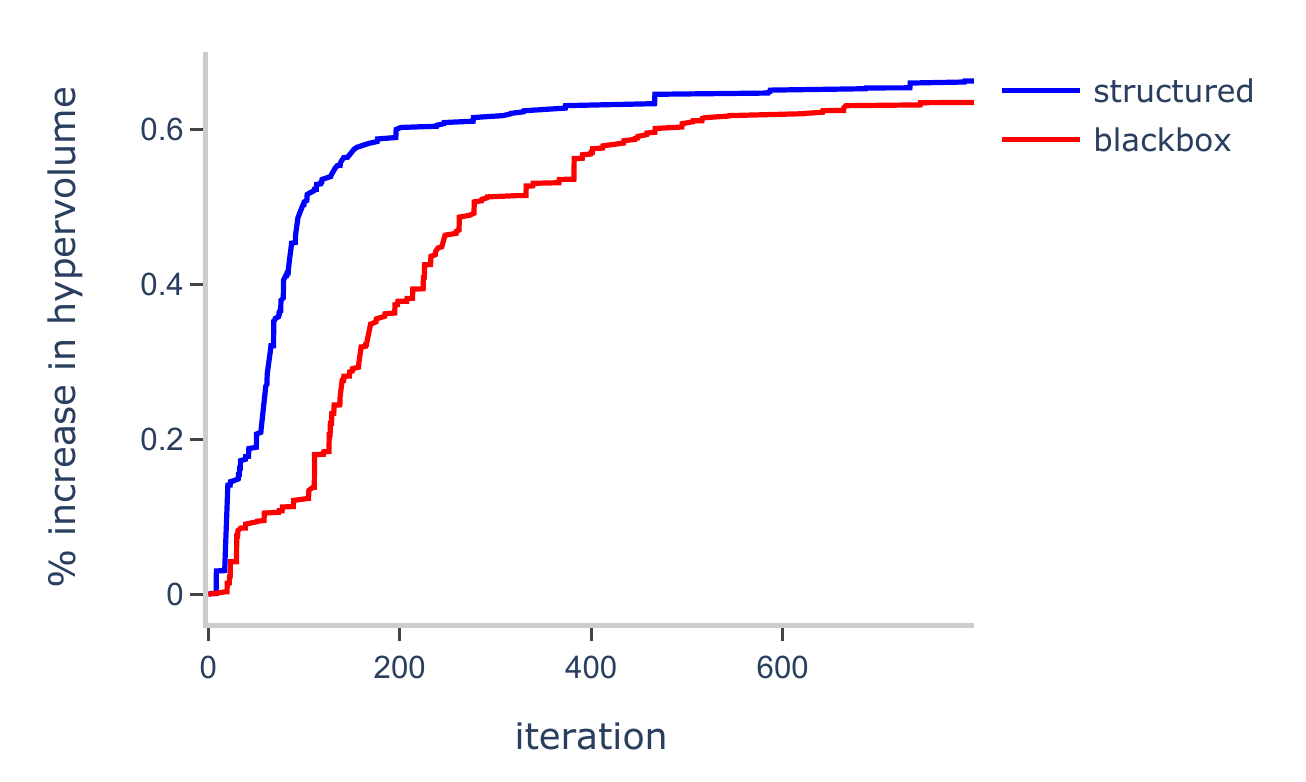}
\end{center}
\caption{Iteration vs.\ percentage hypervolume improvement when solving the
Fayans EDF calibration with ParMOO, exploiting the sum-of-squares
structure (structured) and with a standard (black-box) approach.  The total
simulations used by the end of iteration $k$ are calculated as $2000 + 10k$.
}
\label{fig:fayans-results-hv}
\end{figure}

The results of solving \Cref{eq:fayans-multi} based on the {\tt keras}/{\tt
torch} residual model are presented here.  In order to account for variability
resulting from randomness in the initial design-of-experiments, all performance
results have been averaged across five random seeds.  First, borrowing the
metric used in \cite{bollapragada2020}, we present the $\chi^2$ loss across all
observable classes in \Cref{fig:fayans-results-chi2}.  Note that for $\chi^2$
loss, small values are better.  Next,  to estimate the multiobjective
performance of our methods, we present a rescaling of the hypervolume
performance indicator in Figure \ref{fig:fayans-results-hv}.  For the
hypervolume indicator, large values are better.

Note that the hypervolume indicator is extremely sensitive to problem scaling
and the choice of Nadir point, so its raw values are difficult to interpret.
In an effort to normalize values, we present the improvement in hypervolume
over that of the initial 2,000-point LHS design, as a percentage of the
hypervolume dominated by that original design.  Even after this normalization,
however, the absolute value of the hypervolume improvement is still difficult
to interpret since it is still influenced by the total hypervolume between the
true Pareto front and Nadir point.  In \Cref{fig:fayans-results-hv} we see less
than a 0.7\% increase in total hypervolume.  This is because our Nadir point is
determined by the lower-bound constraints that ParMOO enforced on the range of
interesting values, which were intentionally set to be overly pessimistic.
However, the hypervolume improvements of the two methods relative to each other
can still be taken as an indicator of relative performance.

For both the $\chi^2$ loss and hypervolume improvement, the structured solver
converges considerably faster than does the black-box solver.  With respect to
the $\chi^2$ loss, in just 200 iterations (4,000 simulation evaluations) the
structured solver achieves better performance than the black-box solver will in
all 800 iterations (10,000 simulation evaluations).  This result is to be
expected: since one of our acquisition functions specifically targets this
solution with fixed scalarization weights, ParMOO behaves similarly to how a
structure-exploiting single-objective solver would for this performance metric.

With respect to the percent hypervolume improvement (a true multiobjective
performance metric), the ParMOO's structured solver also achieves significantly
improved performance.  However, the convergence appears to be slightly slower
than for the $\chi^2$ loss.  This is to be expected since solving the full
multiobjective problem is considerably harder than solving a single
scalarization.  However, especially for limited computational budgets, the
performance of the structured solver is still dramatically improved.

\section{Case Study: A Multiobjective Chemical Design}
\label{sec:chem}
To demonstrate ParMOO's effectiveness with heterogeneous objectives and
flexibility to support diverse scientific workflows, in this section we apply
ParMOO to solve a materials engineering problem.  Here, two of the three
objectives are derived from a continuous-flow chemistry experiment that is
conducted in a wet-lab environment, while the third objective is algebraic.  As
in \Cref{sec:fayans}, to make our results open and reproducible, in this
section we use a nonlinear model of the chemical response surface, based on
real-world data collected using nuclear magnetic resonance (NMR) spectroscopy
on the solutions from a continuous-flow reactor (CFR).

\subsection{Background on CFR Material Design}

In material manufacturing applications, our goal is to propose a technique for
chemical manufacturing that can be used to produce a desired material with high
purity at scale.  In particular, in this application we are manufacturing the
electrolyte 2,2,2-trifluoroethyl methyl carbonate (TFMC) by mixing one of two
predetermined potential solvents with one of two predetermined potential bases.
Note that for confidentiality reasons, we do not include the true names of
these potential solvents and bases in this section; they are simply labeled as
S1, S2 and B1, B2, respectively.  This does not affect the reproducibility of
our results since they are also labeled as such in the nonlinear model that we
used to represent the problem in this section.

Our goal is to find an optimal pairing of solvents and bases and conditions
(such as reaction time, equivalence ratio, and temperature) for producing a
pure solution of TFMC in a CFR.  However, in addition to producing a pure
solution, we want to be able to use short reaction times, so that we can
produce large quantities of TFMC at scale.  We expect that this will require us
to use higher temperatures, which could activate a side reaction and produce an
unwanted byproduct, thereby reducing the purity.

\subsection{The Chemical Response Surface Model}
\label{sec:cfr-rs}

In the context of ParMOO, this problem has three continuous design variables
and two categorical design variables.  These variable names and their
ranges/potential values are given in \Cref{tab:cfr-bounds}.  The response
values of interest are the integrals of TFMC and byproduct production values
over a fixed-length time window, as measured by using NMR spectroscopy.

\begin{table}[h]
\begin{center}
\begin{tabular}{c|cccc}
Variable Name   & Var. Type & Lower Bound & Upper Bound & Legal Vals. \\
\hline
temperature (T) & continuous    & 35 C        & 150 C       & N/A          \\
reaction time (RT) & continuous & 60 sec      & 300 sec     & N/A          \\
equiv.~ratio (EQR) & continuous & 0.8    & 1.5         & N/A          \\
solvent         & categorical   & N/A         & N/A         & S1,S2        \\
base            & categorical   & N/A         & N/A         & B1,B2        \\
\end{tabular}
\vskip 12pt
\end{center}
\caption{ Variable types, bounds, and legal values for the CFR material
optimization problem discussed in this section.  Note that continuous variables
have bounds, while categorical variables have legal values.  }
\label{tab:cfr-bounds}
\end{table}

Since we cannot provide access to the CFR for our experiments, in this section
we use a nonlinear response surface representing the physical continuous-flow
chemistry experiment.  These models were fit by using real-world data that was
collected by providing a physical CFR/NMR setup as a simulation that ParMOO
could query in closed loop.  To do so, ParMOO's extensible API was
layered over the MDML tool \cite{elias2020}, which uses an Apache kafka backend
to distribute requests for experiment evaluations to the CFR and collects and
returns NMR results directly to an online database that ParMOO can query for
simulation results.  For more information on how this data was collected, see
\cite{chang2023c}, which describes the collection of a smaller dataset by using
an identical experimental setup.

After a budget of 62 experiments, ParMOO had converged on several approximate
solutions to the real-world problem, and it was no longer economically viable
to continue the real-world experiment since the cost of real-world materials is
nonnegligible.  The resulting experimental database of design
point/integral-value pairs was used in this section to fit the nonlinear models
described above.  Since this dataset is relatively sparse in the 5-dimensional
input space, special care was taken to ensure that the resulting model does not
exhibit non-physical behaviors.  In particular, to verify our model, we have
ensured that both of our response surfaces
\begin{itemize}
\item approximate the underlying data with low MAE, particularly
for near-optimal design points;
\item do not take on negative values anywhere in the feasible design space,
which would be physically impossible since the outputs represent time integrals
of material production;
\item do not have a sum that exceeds the total amount of solvent and base
provided for any inputs in the design space; and
\item do not take unexpected maxima/minima along the boundaries of the design
space, which could be an artifact of a lack of data in those regions, allowing
our model to overextrapolate.
\end{itemize}

To achieve these criteria, we hand-crafted a small number of physically
meaningful nonlinear terms for our response surface based on the expected
chemistry of the reaction.  We then fit the coefficients of these terms to each
of the two integrals, by using a combination of generalized linear regression
with {\tt scikit-learn} \cite{pedregosa2011}.  The resulting MAE was found to
be within the acceptable ranges.  Using hyperparameter tuning with the {\tt
Powell} solver from {\tt scipy.optimize.minimize} \cite{virtanen2020}, we were
also able to guarantee that the individual global minima for each model was
nonnegative and the sum of the global maxima was within the acceptable range.
We then verified that the individual minima and maxima for each model were
located in acceptable regions of the design space, which agree with our
physical intuition and empirical data.

In the remainder of this section we use these trained response surface models
as the true chemical response surfaces for both the TFMC and byproduct
integrals.

\subsection{Heterogeneous Problem Structure}

In this problems, two of the three objectives are the result of a true
black-box experiment, which must be carried out in a wet-lab environment using
real materials and with  an extremely restrictive budget.  However, the third
objective represents the reaction time, which is one of the directly
controllable inputs.  Therefore, ParMOO is able to directly control this
output, and is able to exploit
this ability to accelerate its practical convergence.

\subsection{The ParMOO CFR Material Design Solver}

To define the problem, we provide ParMOO with the five design variables shown
from \Cref{tab:cfr-bounds}, with their respective variable types and
bounds/values.  To handle the two categorical variables, we use ParMOO's
default categorical variable embedder, which embeds the four distinct
combinations of categorical variables into a three-dimensional continuous
latent space.  Combined with the three continuous design variables, this
results in a six-dimensional effective optimization space.

For the structured variation of the problem, ParMOO is provided with the
pretrained nonlinear model of the chemical response surface,  described in
\Cref{sec:cfr-rs}.  ParMOO is configured to treat this as a single simulation
with two outputs, using a Latin hypercube sample (LHS) with 50 evaluations to
produce the initial database and Gaussian RBFs for surrogate modeling.  Next,
ParMOO is given three objectives, two of which are identity mappings from the
simulation outputs and the third of which is an identity mapping from the
reaction time design variable.  Then, two epsilon-constraint 
functions and one fixed-weighting are added, resulting in a batch size of three
simulation evaluations per iteration.

For comparison and to demonstrate the advantage in exploiting the heterogeneous
problem structure, we also provide an identical implementation of ParMOO for
this problem, where the third ``reaction time'' objective is provided to ParMOO
as a third black-box simulation output.  This results in ParMOO modeling the
third output as a black-box function and ignoring the heterogeneous structure,
as would occur when using most other off-the-shelf multiobjective black-box
optimization solvers.  All other settings are identical as in the structured
variation defined above.

One of the critical challenges of performing automatic experimentation in the
context of material manufacturing is the cost of raw materials and performing
real-world experiments.  Although we do not have these costs for our
computational model of the material response surface, we have run ParMOO with a
restrictive budget of just a $50$-point initial design followed by $30$
iterations with $3$ acquisition functions ($140$ total experiments).  By
exploiting the heterogeneous problem structure, we hope to still achieve good
performance (especially for our cheap objective) with this limited budget.

The response surface model described above and the code for reproducing the
experiments presented here are given at
\url{https://github.com/parmoo/parmoo-solver-farm/tree/main/cfr-material-design-2022}.
The results are presented in the next section.

\subsection{Results}
\label{sec:cfr-results}

The results of tuning manufacturing conditions based on our chemical response
surface are presented here.  Again, all performance results have been averaged
across 5 unique random seeds.

\begin{figure}[h!]
\begin{center}
\includegraphics[width=0.75\textwidth]{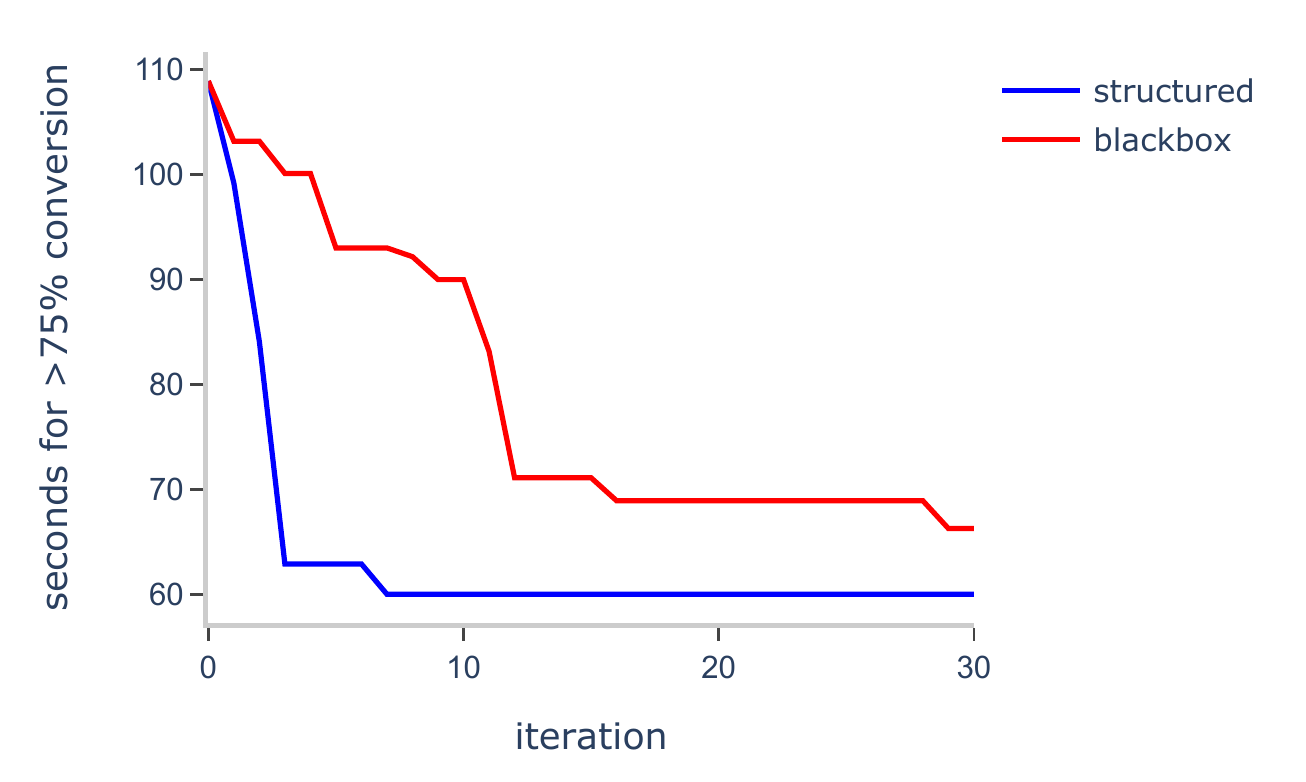}
\end{center}
\caption{Iteration vs.\ minimum reaction time that achieves at least a 75\%
conversion rate when solving the CFR chemical manufacturing problem with
ParMOO, exploiting the heterogeneous structure (structured) and with a
standard (black-box) approach.  The total simulations used by the end of
iteration $k$ are calculated as $50 + 3k$.  }
\label{fig:cfr-results-rt}
\end{figure}

\begin{figure}[h!]
\begin{center}
\includegraphics[width=0.75\textwidth]{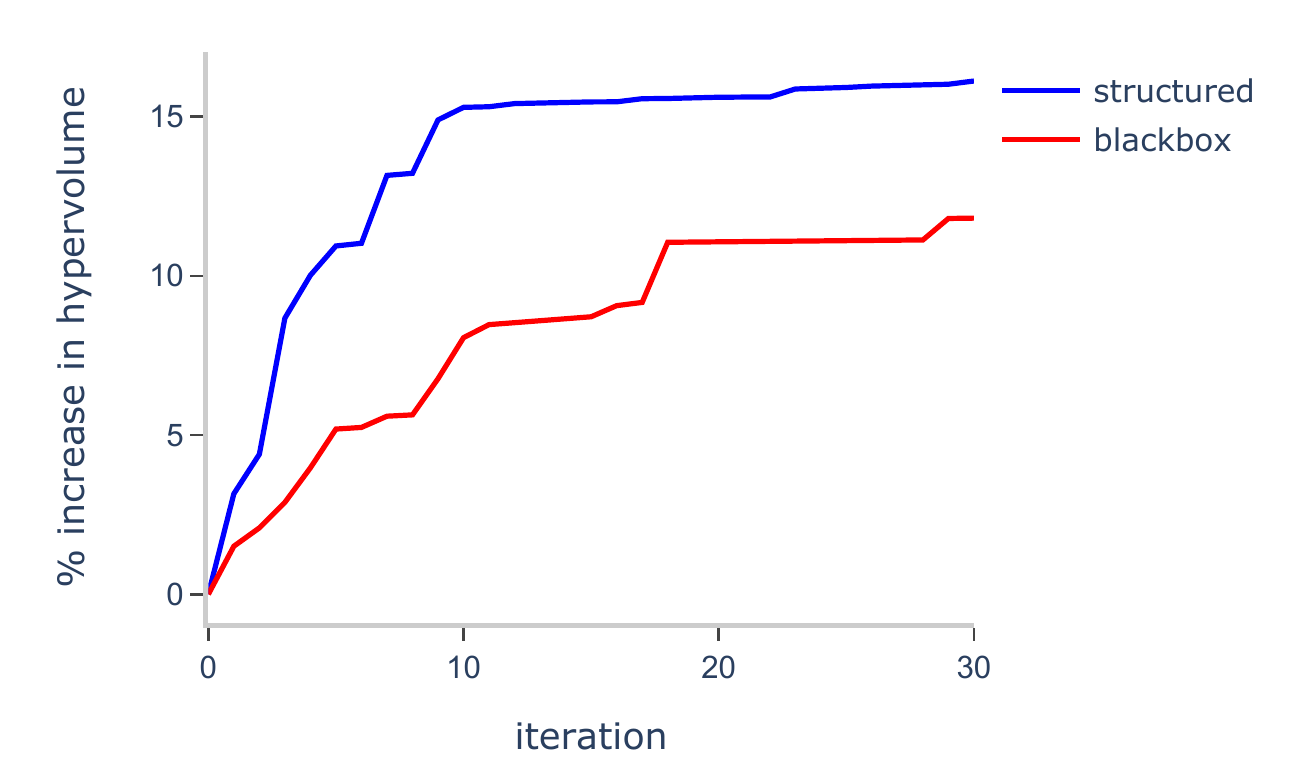}
\end{center}
\caption{Iteration vs.\ percentage hypervolume improvement when solving the CFR
chemical manufacturing problem with ParMOO, exploiting the heterogeneous
structure (structured) and with a standard (black-box) approach.  The total
simulations used by the end of iteration $k$ are calculated as $50 + 3k$.  }
\label{fig:cfr-results-hv}
\end{figure}

First,  to assess our performance with respect to our one cheap objective, we
present the average minimum observed reaction times (in seconds) subject to a
75\% purity constraint in \Cref{fig:cfr-results-rt}.  Note that for these
reaction times, small values are better.  Next,  as in
\Cref{sec:fayans-results}, we present the improvement in hypervolume as a
proportion of the gap between initial hypervolume and total possible
hypervolume (with respect to the ideal point) in \Cref{fig:cfr-results-hv}.
Again, for the hypervolume improvement, large values are better.

As in \Cref{sec:fayans-results}, the structured solver greatly outperforms the
black-box approach by both metrics.  We note that although the time required to
achieve over 75\% conversion is a physically meaningful convergence metric and
demonstrative of ParMOO's ability to directly control the reaction time (the
algebraic objective), this time we did not provide any fixed scalarization that
would explicitly target this solution, reducing to a single-objective problem.

We note that ParMOO's structure-exploiting solver achieves excellent
performance on both this problem and the Fayans problem from
\Cref{sec:fayans-results} with little additional work from the user, even
though these structures are considerably different.

\section{Discussion and Continued Work}
\label{sec:discussion}

In this paper we have described the design principles behind the design of the
MOSO library ParMOO.  To summarize, our five main design goals are
\begin{enumerate}
\item customizability of solvers;
\item exploitation of
composite
structures in MOSO problem formulation;
\item flexibility in support for a wide variety of design spaces;
\item ease of use in scientific workflows;
and
\item usability, extensibility, and maintainability as an open-source software
package.
\end{enumerate}

We have achieved these goals through an object-oriented design that is highly
modularized, utilizes an intermediate simulation output space, and embeds
complex problems into a continuous latent input space.  This framework has been
demonstrated on open-source models of two real-world problems.
These models and the code for reproducing our results have been shared through
our {\tt parmoo-solver-farm} repository, and a snapshot of all related code at
the time of publication is available in \cite{chang2024}.

Although these two problems exhibit completely different structures, ParMOO is
able to effectively exploit the structure in both cases and  to considerably
improve the  convergence relative to that of a black-box approach that does not
exploit known structure.  Notably, defining such a structure-exploiting solver
requires little additional work or customization from the user.

\section*{ACKNOWLEDGMENTS}
We are grateful to Jared O'Neal, Witold Nazarewicz, and Paul-Gerhardt Reinhard
for providing the Fayans functional data and to Joseph Libera, Jakob Elias,
Santanu Chaudhuri, and Trevor Dzwiniel for their respective roles in the
collection
of the
CFR materials data.  We are also grateful to Jeff Larson, John-Luke Navarro,
and Steve Hudson for their advice on best practices in scientific software
development and support in our usage of libEnsemble.  We are grateful to four
anonymous reviewers for their comments, which have improved the presentation of
this work.

This work was supported in part by the U.S.~Department of Energy, Office of Science, Office of Advanced Scientific Computing Research's SciDAC program under Contract Nos.\ DE-AC02-05CH11231 and~DE-AC02-06CH11357.

\bigskip
\bibliographystyle{elsarticle-num-names}
\bibliography{moo-refs.bib}

\begin{flushright}
\scriptsize
\framebox{\parbox{\textwidth}{
The submitted manuscript has been created by UChicago Argonne, LLC, Operator of
Argonne National Laboratory (“Argonne”).  Argonne, a U.S. Department of Energy
Office of Science laboratory, is operated under Contract No. DE-AC02-06CH11357.
The U.S. Government retains for itself, and others acting on its behalf, a
paid-up nonexclusive, irrevocable worldwide license in said article to
reproduce, prepare derivative works, distribute copies to the public, and
perform publicly and display publicly, by or on behalf of the Government.  The
Department of Energy will provide public access to these results of federally
sponsored research in accordance with the DOE Public Access Plan.
\url{http://energy.gov/downloads/doe-public-access-plan}
}}
\normalsize
\end{flushright}

\end{document}